\documentclass[journal ]{new-aiaa}
\usepackage[utf8]{inputenc}
\usepackage{textcomp}
\usepackage{todonotes}
\usepackage{graphicx}
\usepackage{amsmath}
\usepackage[version=4]{mhchem}
\usepackage{siunitx}
\usepackage{longtable,tabularx}
\usepackage[linesnumbered,ruled,vlined]{algorithm2e}
\usepackage{rotating}
\usepackage{pdflscape}
\usepackage{multirow}

\SetCommentSty{mycommfont}
\title{Hierarchical Framework for Space Exploration Campaign Schedule Optimization\footnote{A previous version of this paper was presented at the AIAA SciTech Forum in National Habor, MD \& Virtual on January 23-27, 2023 (AIAA 2023-1965).}}

\author{Nicholas Gollins\footnote{PhD Candidate, Daniel Guggenheim School of Aerospace Engineering,  AIAA Student Member.} and Koki Ho\footnote{Dutton-Ducoffe Professor, Associate Professor, Daniel Guggenheim School of Aerospace Engineering, AIAA Senior Member. kokiho@gatech.edu}}
\affil{Georgia Institute of Technology, Atlanta, Georgia, 30332}

\begin{document}

\maketitle

\begin{abstract}
Space exploration plans are becoming increasingly complex as public agencies and private companies target deep-space locations, such as cislunar space and beyond, which require long-duration missions and many supporting systems and payloads. Optimizing multi-mission exploration campaigns is challenging due to the large number of required launches as well as their sequencing and compatibility requirements, making the conventional space logistics formulations not scalable. To tackle this challenge, this paper proposes an alternative approach that leverages a two-level hierarchical optimization algorithm: a genetic algorithm is used to explore the campaign scheduling solution space, and each of the solutions is then evaluated using a time-expanded multi-commodity flow mixed-integer linear program. A number of case studies, focusing on the Artemis lunar exploration program, demonstrate how the method can be used to analyze potential campaign architectures. The method enables a potential mission planner to study the sensitivity of a campaign to program-level parameters such as logistics vehicle availability and performance, payload launch windows, and in-situ resource utilization infrastructure efficiency. 
\end{abstract}

\section*{Nomenclature}

{\renewcommand\arraystretch{1.0}
\noindent\begin{longtable*}{@{}l @{\quad=\quad} l@{}}
$c$                 & Crew consumables consumption rate (kg per crew member per day) \\
$C$                 & Mixed-integer linear program objective function cost multiplier \\
$\mathcal{C}$       & Set of co-payloads \\
$d$                 & Demand matrix \\
$\mathcal{D}$       & Spacecraft domain \\
$\mathcal{E}$       & Boolean variable defining whether an arc exists. \\
$f$                 & Objective function of commodity flow mixed-integer linear program \\
$\mathcal{F}$       & Objective function of metaheuristic optimizer and output value of mixed-integer linear program commodity flow function \\
$i$                 & Origin node \\
$I$                 & Total number of logistics network nodes \\
$io$                & Flow direction (in-to- or out-of-arc) \\
$I_{\mathrm{sp}}$            & Specific impulse (s)\\
$j$                 & Target node \\
{$k$}    & {Penalty function scaling coefficient} \\
$m$                 & Quantity or mass (kg)\\
$m_{\mathrm{dry}}$           & Dry Mass (kg)\\
{$m_{\mathrm{crew}}$}           & {Mass per Crew Member (kg per person)}\\
$m_{\mathrm{pay}}$           & Payload capacity (kg)\\
$m_{\mathrm{prop}}$          & Propellant capacity (kg)\\
$n$                 & vehicle index \\
$N$                 & Total number of vehicle designs \\
${\mathcal{N}}$  & {Number of duplicated time steps in the nonlinear penalty case study} \\
$p$                 & Payload type index \\
$P$                 & Total number of commodity types \\
$P_F$               & Number of float (continuous) commodity types \\
$P_I$               & Number of integer commodity types \\
$\mathcal{P}$       & Set of soft precursor payloads  \\
$\mathcal{Q}$       & Set of strict precursors payloads \\
$\mathcal{R}$       & Set of programmatic requirements \\
$\mathcal{S}$       & Set of vehicles that can form a stack \\
$\mathcal{S}'$      & Set of stacks to which a vehicle belongs \\
$t$                 & Time index \\
$t_F$               & Minimum time between launches (months)\\
$t_L$               & Lower bound of  launch window \\
$t_{U}$             & Upper bound of launch window \\
$T_{\mathrm{LP}}$            & Total number of time steps in the linear program \\  
$x$                 & Commodity flow (decision variable of mixed-integer linear program function) \\
$\textbf{x}$        & Metaheuristic decision vector and input variable to $\mathcal{F}$ \\
$Z$                 & Propellant mass fraction \\
$\beta$             & Propellant boil-off rate  ($\%$ per day)\\
$\Gamma$ & Metaheuristic cost function \\ 
$\Lambda$        & Metaheuristic penalty function \\
$\phi$              & Propellant oxidizer mass ratio \\
$\mu$               & in-situ resource utilization maintenance supply requirement (kg per kg ISRU infrastructure per day) \\
$\rho$              & in-situ resource utilization propellant production rate (kg propellant per kg in-situ resource utilization infrastructure per day) \\
$\tau$              & Time of flight \\
\end{longtable*}}

\section{Introduction}

\lettrine{T}{he} Global Exploration Roadmap (GER), published by the International Space Exploration Coordination Group (ISECG) \cite{GER}, lays out the development of a sustained crewed presence in beyond-low Earth orbit (LEO) locations over the coming decades. The document focuses particularly on lunar surface exploration and the planned Artemis program and supporting missions. Ref. \cite{GER} shows a possible sequence of missions in the Artemis and supporting programs, demonstrating the increasing cadence of missions launched to cislunar space over the coming decade. Historically, however, the crewed exploration of space has taken one of two forms:

\begin{itemize}
    \item Early LEO missions, through the Apollo era, and the early Shuttle era, largely consisted of single-shot missions that carried all of the required material in a single launch.
    \item The International Space Station (ISS) era required a more sophisticated logistics plan, with a variety of resources sourced from different providers being delivered to the space station.
\end{itemize}

The plans laid out in the GER combine the challenges of deep-space exploration of the Apollo era with the complex logistics requirements of the ISS. This combination of increasing complexity and cadence of missions, many of which are interdependent, manufactures a necessity for efficient planning. It would be useful, therefore, for a would-be campaign architect to be able to find an optimal campaign plan for a particular set of programmatic and system assumptions, and to be able to study how the optimal plan reacts to changes in those assumptions. The optimization method presented in this paper is intended to serve as a process by which to obtain a well-performing baseline campaign plan according to pre-determined requirements.

There have been various studies in the space logistics research area to tackle such challenges. Network-based optimization, along with mixed-integer linear programming (MILP), was used to model ISS resource and traffic management and used to assess the usefulness of logistics vehicle designs \cite{Fasano1996, Fasano1997}. More recently, the field has expanded to address the problems of modern space systems optimization problems. Taylor \textit{et al} \cite{Taylor_Klabjanz_Simchi-levi_Song_DeWeck_2006} developed a MILP-based heuristic commodity flow model for interplanetary transport, in which commodity flow paths are pre-determined and then assigned to logistics vehicles. Later, the SpaceNet space logistics modeling framework utilized time-expanded networks and was used to study ISS, lunar, and interplanetary logistics \cite{Gralla2006, Grogan_Yue_DeWeck_2011}. More recently, multi-commodity flow models have been developed to calculate the optimal flow of spacecraft, payloads, and propellant through space \cite{Ho_deWeck_Hoffman_Shishko_2014, Ishimatsu2016}. Furthermore, space logistics models have been developed further to incorporate spacecraft design into the overall mission optimization framework \cite{Chen_Ho_2018, Isaji2021, Chen2021, Takubo2022}. 
In practice, uncertainty can arise in many aspects of space mission planning. For example, Blossey \cite{Blossey_2023} constructed a stochastic MILP to tackle commodity demand uncertainty in Mars mission planning. Other examples of sources of uncertainty could result from launch delay, in-space maneuver error, operational risks, and/or system failures. Chen \textit{et al} established a method for adjusting an optimal commodity flow, for a given schedule, for robustness to launch delays \cite{Chen_Gardner_Grogan_Ho_2021}. 

While most of the existing space logistics mission design studies assumed a relatively well-defined set of requirements and hypothetical/simplified vehicles, in practice, space mission designers need to explore a large design space with a set of requirements that are not necessarily precisely defined and yet interdependent. For example, a multi-mission campaign needs to be scheduled with large possible launch windows while satisfying the payload delivery sequencing requirements and availability of each potential vehicle. Considering the large number of payloads that can feature in multi-mission exploration campaigns, there can be many combinations of choices of launch times for each payload within their respective launch windows, each with many choices of available vehicles. Different combinations can result in different logistics solutions, some more efficient than others. For smaller campaigns, an obvious solution might be, for example, to group all payloads into a single, larger spacecraft as a ``ride-share'' mission. However, the more payloads and possible vehicles that are involved in the campaign, the more difficult it becomes to navigate the solution space efficiently. This is due to the complex combinations of programmatic requirements and vehicle availability constraints creating a sparse solution space. It would be beneficial, therefore, to have a method by which to efficiently assess commodity flow parameters and model structures, such as launch schedules or vehicle properties. 

While some of these complexity challenges have been tackled by MILP in conventional studies, these methods have significant limitations in terms of their scalability to the mission sequencing requirements and the number/types of payloads/vehicles. For example, the lower and upper bounds of payload launch times and vehicle availability can be translated into supply/demand times within the commodity flow in the conventional MILP formulation. However, more complex types of constraints such as payload or mission sequencing (e.g. payload X must launch before payload Y) place constraints on the interactions between the commodities within the MILP. This can be handled in the MILP by treating each individual payload as a separately-indexed commodity type, and tailoring sequencing constraints to each of those commodities, but creating more commodity types increases the size of the problem, and solving times can grow rapidly or the solver can encounter memory limitations. \footnote{The scalability issues of MILP are discussed further in Appendix A.} While other approaches exist by introducing auxiliary variables to enforce mission sequencing \cite{sarton2021framework}, those approaches still increase the numbers of (integer) variables and constraints and thus complexity. In addition to this, the solution space is sparse: a large number of infeasible solutions can lie between two feasible ones, which poses challenges to the solver. 

Another drawback of using MILP is that it cannot capture a nonlinear cost function in commodity flow modeling. Methods exist for solving commodity flow problems with nonlinear costs appearing in the objective function \cite{Goffin_Gondzio_Sarkissian_Vial_1996} and/or with nonlinear constraints \cite{Castro_Nabona_1996, Mijangos_2006}. However, these methods only work with continuous variables - mixed-integer nonlinear commodity flows are much more difficult to solve, but it is possible using metaheuristic algorithms \cite{Alavidoost_Tarimoradi_Zarandi_2018}. Mixed-integer formulations are necessary in space logistics models due to the need to model the transport of discrete commodities such as spacecraft and crew.

To tackle the above challenges, this work aims to propose a generalized schedule optimization by using a multi-level optimization. Namely, a genetic optimization algorithm is used to explore scheduling solutions, which are evaluated using a logistics MILP formulation. Metaheuristic optimization has been utilized in previous space logistics works. Taylor \textit{et al} used a heuristic random search algorithm to navigate solution to a commodity flow problem \cite{Taylor_Klabjanz_Simchi-levi_Song_DeWeck_2006}; Chen and Ho \cite{Chen_Ho_2018} used simulated annealing optimization to integrate vehicle design with commodity flow; and Jagannatha and Ho \cite{Jagannatha_Ho_2018} used a multi-objective genetic algorithm to find optimal parameters for in-space propellant depot supply chain design. In this work, the intended use of the genetic scheduling algorithm is to schedule the demand times of program payloads such as crew, habitats, rovers, or scientific equipment, for example, considering nonlinear constraints. At the same time, the MILP handles the "supporting" payloads such as crew supplies and maintenance equipment. This limits the number of payload-type indices required and therefore keeps the scale of the MILP problem manageable. An advantage of using a metaheuristic outer level is that the overall objective can become a nonlinear cost function as long as the nonlinearity is caused by the decision variables of the metaheuristics level. 
This capability allows a much broader range of cost functions; for example, it would allow program planners to apply a nonlinear penalty to the MILP solution as a function of the launch schedule, which would not be possible in a MILP-only approach.  
Another advantage of the proposed formulation is that it could be structured to facilitate uncertainties in the decision vector due to its generality in the form of the cost functions, though this is left for future work.

The results are validated and demonstrated with realistic case studies based on NASA's Apollo, Commercial Lunar Payload Services (CLPS) programs, and Artemis programs. These programs were chosen as case studies because of their complexity in terms of the wide range of payloads and logistics service providers. With NASA's recent approach of contracting out payload delivery to logistics service providers, it is important to maintain a centralized perspective on the available logistics vehicles in order to ensure that payload delivery contracts are awarded according to an efficient plan. This is particularly true when many of those payloads are in support of NASA's crewed exploration program, for which NASA maintains operational responsibility.

The rest of the paper is organized in the following way. First, Section \ref{sec:Method} describes the two-layer optimization algorithm in detail, with Section \ref{sec:Meta} focusing on the metaheuristic schedule optimizer and \ref{sec:MILP} focusing on the commodity flow MILP. Then, the method is tested on a number of case studies in Section \ref{sec:Results}. The commodity flow MILP is demonstrated in Section \ref{sec:Apollo} using programmatic requirements based on an Apollo mission as input. In Section \ref{sec:CLPS}, NASA's CLPS program is studied as a demonstration of the scheduling algorithm. Finally in Sections \ref{sec:Art1-2A} and \ref{sec:Art2B}, the first phases (1, 2A, and 2B) of the Artemis program according to the GER are studied.

\section{Method}
\label{sec:Method}
\subsection{Overview}
The overarching method of the campaign planner is as follows: 

\begin{itemize}
    \item Vehicle design and campaign programmatic requirements are defined in input databases. The structure of the data is defined in Sections \ref{sec:Vehicle_Data} and \ref{sec:Program_Data}. 
    \item A series of initial feasible guesses of the scheduling solutions are generated and used to initialize a genetic algorithm population. The decision variable contains the launch time of each payload in the campaign and this determines the network and the ``demand matrix'' of the logistics model.
    \item  With each solution assessment, a MILP is constructed to solve the logistics optimization problem. If the model is feasible, it returns the total launch mass across the campaign as the objective to the genetic algorithm.
\end{itemize}

\begin{figure}[]
\centering
\includegraphics[width=1\textwidth]{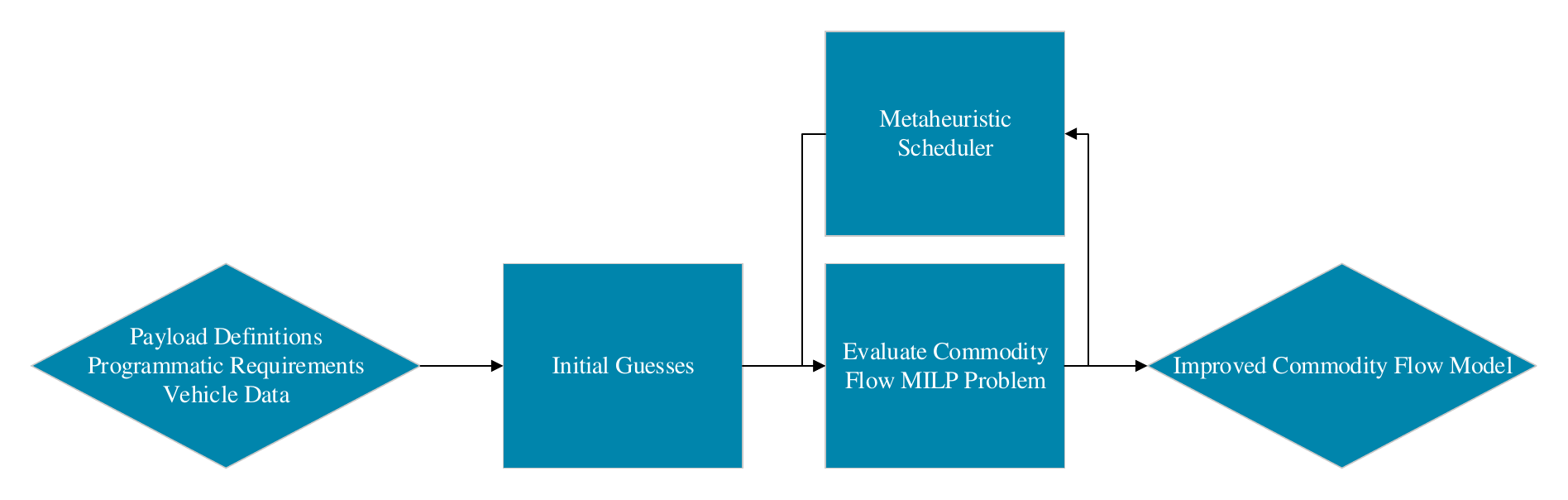}
\caption{Overall structure of the generic campaign optimization problem \cite{Gollins_Isaji_Shimane_Ho_2023}.}
\label{problem_structure}
\end{figure}

The evolutionary algorithm repeats this process in search of more optimal (minimized total launch mass) solutions. This workflow is illustrated in Figure \ref{problem_structure}.

\subsubsection{Vehicle Data} \label{sec:Vehicle_Data} The vehicle data input database file contains a list of logistics vehicles that are available for the lunar exploration campaign. Table \ref{tab:vehicle_variables} lists the parameters associated with each vehicle. The "secondary" parameters are determined based on the values of $I_{\mathrm{sp}}$, which is used as a proxy for propellant type. $I_{\mathrm{sp}} < 370$ was taken to be a storable propellant type, $I_{\mathrm{sp}}\geq370$ s was taken to be cryogenic propellant. In this paper, our focus is on the chemical propulsion system, but the work can be extended to other propulsion systems such as electric propulsion and nuclear thermal propulsion systems\footnote{Note that for low-thrust electric propulsion, the trajectory design needs to be incorporated into the logistics design \cite{Jagannatha_Ho_2020}.}.  Throughout all case studies, the oxidizer boil-off rate was 0.016$\%$ loss per day for liquid oxygen and 0 for storable propellants. Fuel boil-off was neglected. Table \ref{tab:example_vehicles} gives an example of vehicle data input format, in which Vehicle 1 is a small lander and Vehicle 2 is a larger service module.

In addition to the list of vehicles, a list of vehicle ``stacks'' is defined. A ``stack'' is an aggregation of individual space vehicles into a single unit, that can be assembled through rendezvous or disconnected. The vehicle stacking system implemented in the linear program was inspired by the ontologies developed by Trent, Edwards \textit{et al}, and Downs \cite{Trent_2017, Edwards_Trent_Diaz_Mavris_2018, Downs_Prasad_Robertson_Mavris_2023}. A vehicle stack is formatted as a list of its constituent vehicles and creates a new vehicle with the $I_{\mathrm{sp}}$ and domain of the leading element, which is taken to be the active element. For example, referring to the example vehicles in Table \ref{tab:example_vehicles}, if the stack [2,1] is allowed, then a new Vehicle 3 becomes the stack and has the domain and $I_{\mathrm{sp}}$ of Vehicle 2. Dry masses are summed together. The payload capability of the stack is the payload capability of the leading element minus the dry masses of the inactive elements. 

\begin{table}[]
	\fontsize{10}{10}\selectfont
    \caption{Definitions of the input vehicle parameters.}
   \label{tab:vehicle_variables}
        \centering 
   \begin{tabular}{c l } % Column formatting, 
      \hline 
      Parameter    & Description \\
      \hline 
      $m_{\mathrm{pay}}$   & Payload capacity \\
      $m_{\mathrm{prop}}$  & Propellant capacity \\
      $m_{\mathrm{dry}}$   & Dry Mass \\
      $I_{\mathrm{sp}}$    & Specific impulse \\
      $t_F$     & Launch frequency \\
      $t_L$     & Earliest allowed launch \\
      $\mathcal{D}$     & Domain: the set of arcs along which the vehicle is allowed to travel. \\
      \hline 
      Secondary Parameter    & Description \\
      \hline 
      $\phi$    & Propellant oxidizer mass ratio \\
      $\beta$   & Boil-off rate  \\
      \hline
   \end{tabular}
\end{table}

\begin{table}[]
    \centering
    \caption{Example vehicle data input.}
    \begin{tabular}{l c c c c c c c}
        \hline
         Name & $m_{\mathrm{pay}}$, kg & $m_{\mathrm{prop}}$, kg & $m_{\mathrm{dry}}$, kg & $I_{\mathrm{sp}}$, s & $t_F$ & $t_L$ & $\mathcal{D}$ \\
         \hline
         Vehicle 1 & 90 & 720 & 470 & 340 & 12 & 1 & [0,0], [0,1], [1,2], [2,2], [2,3], [3,3] \\
         Vehicle 2 & 1400 & 3320 & 1950 & 340 & 12 & 24 & [0,0], [0,1], [1,2], [2,2] \\
         \hline
    \end{tabular}
    
    \label{tab:example_vehicles}
\end{table}

\subsubsection{Programmatic Requirements} \label{sec:Program_Data}
The campaign program definition database file contains a list of payloads to be launched throughout the campaign. Table \ref{tab:campaign_variables} lists the parameters associated with each payload. Payloads are indexed by $l \in [0, N_p)$, where $N_p$ is the total number of campaign payloads. So, the first payload is $l=0$, for example. The sets of pre- and co-payloads contain the indices of those payloads. The type indices are the same as those used in \cite{Chen_Ho_2018, Gollins_Isaji_Shimane_Ho_2023}. Table \ref{example_payloads} gives an example of the payload data input format.

\begin{table}[]
	\fontsize{10}{10}\selectfont
    \caption{Definitions of the programmatic parameters.}
   \label{tab:campaign_variables}
        \centering 
   \begin{tabular}{c l } % Column formatting, 
      \hline 
      Parameter    & Description \\
      \hline 
        $p_l$         & Payload type \\
        $m_l$       & Quantity/mass \\
        $i_l$         & Start node \\
        $j_l$         & Target node \\
        $t_{L,l}$     & Lower bound of launch window \\
        $t_{U,l}$     & Upper bound of launch window \\
        $\mathcal{P}_l$     & Set of soft precursors (payload $l$ must arrive \textit{after or with} payloads in this set) \\
        $\mathcal{Q}_l$     & Set of strict precursors (payload $l$ must arrive \textit{strictly after} payloads in this set) \\
        $\mathcal{C}_l$     & Set of co-payloads (payload $l$ must launch \textit{with} payloads in this set) \\
      \hline

   \end{tabular}
\end{table}

\begin{table}[]
    \centering
    \caption{Example payload data input.}
    \begin{tabular}{l c c c c c c c c c c}
        \hline 
         $\#$ & Name & Type Index & Quantity & Supply Node $i$ & Demand Node $j$ & $t_L$ & $t_U$ & Co-Payloads \\
         \hline
         0 & Crew & 1 & 2 & 0 & 3 & 0 & 12 & \\
         1 & Science Payload & 5 & 14 & 0 & 3 & 0 & 6 & 0 \\
         \hline
    \end{tabular}
    \label{example_payloads}
\end{table}

\subsection{Metaheuristics Layer}
\label{sec:Meta}
The ``outer layer'' of the optimization algorithm searches for optimal launch schedules. A genetic algorithm is used, where the decision vector $\textbf{x}\in \boldsymbol{Z}^{N_P}$ contains the launch time of each payload in the campaign. The integer time steps chosen by the metaheuristics are a month in length. This allows for simplified synergy with the periodic nature of the low-energy transfers. \textbf{x} has the form $ \textbf{x} = [t_0, t_1,... t_{N_p}] $

The genetic algorithm searches for solutions to the problem stated in Equation (\ref{eq:meta_obj}), where $\mathcal{F}$ is the output of the MILP, described in Section \ref{sec:MILP}, given the input schedule \textbf{x} and set of programmatic requirements $\mathcal{R}$. $\Lambda$ is a scaling factor on the MILP objective and is a function of the schedule, and $\Gamma$ is a cost associated with the schedule. The use of a metaheuristic optimizer allows both $\Lambda$ and $\Gamma$ to be non-linear if necessary. For simplicity in the case studies, a one-to-one mapping between the MILP and metaheuristic objectives was used, i.e., $\Gamma = 0$ and $\Lambda = 1$. Section \ref{sec:nonlinear} demonstrates a scheduling case with nonlinear $\Gamma(\mathbf{x})$.

\begin{equation} \label{eq:meta_obj}
    \begin{split}
        \min_{t_l} & \quad \Lambda(\textbf{x}(t_l), \mathcal{R})\mathcal{F}(\textbf{x}(t_l), \mathcal{R}) +  \Gamma(\textbf{x}(t_l), \mathcal{R})\\
        \text{s.t.} & \quad t_{L,l} \leq t_l \leq t_{U,l} \\
    \end{split}
\end{equation}

The metaheuristics problem is subject to constraints imposed by the programmatic payload-sequencing requirements. Equations (\ref{eq:meta_constr1}) - (\ref{eq:meta_constr3}) list the constraints corresponding to the soft precursor, strict precursor, and co-payload requirements respectively. Soft precursors $q$ must launch after, or with, the payload $l$, whilst strict pre-cursors must launch strictly after payload $l$. Co-payloads must launch at the same time as payload $l$.

\begin{equation} \label{eq:meta_constr1}
    t_q - t_l \leq 0 \quad \forall \; q \in \mathcal{P}_l       
\end{equation}
    
\begin{equation}\label{eq:meta_constr2}
    t_q - t_l < 0 \quad \forall \; q \in \mathcal{Q}_l
\end{equation}

\begin{equation}\label{eq:meta_constr3}
     t_q - t_l = 0 \quad \forall \; q \in \mathcal{C}_l 
\end{equation}

The constraints are addressed indirectly by applying a ``death penalty'', where a large positive scalar is assigned to the objective function value if the solution is infeasible, regardless of the extent of constraint violation \cite{DeathPenalty}. Similarly, the death penalty is applied to any solution that, although feasible according to the metaheuristics constraints, is found to be infeasible by the MILP layer. The metaheuristic-layer optimization was carried out using the genetic algorithm of the pygmo \cite{Biscani2020} metaheuristic optimization python library.

\subsection{Logistics Mixed-Integer Linear Program (MILP) Layer} \label{sec:MILP}
In space logistics, the transfer of vehicles, crew, payloads, and propellant are modeled as commodities flowing through a network. A mixed-integer linear program, an extension of the models described in \cite{Chen_Ho_2018, Gollins_Isaji_Shimane_Ho_2023}, is used to find the flow that results in a minimized launch mass, subject to certain constraints. Table \ref{tab:model_variables} lists the constants, indices, and variables used in the linear program. Table \ref{tab:nodes} lists the locations represented in the model, and Table \ref{tab:payloads} lists the types of commodities and their costs. The cost of each commodity is its mass per unit. 

\begin{table}[]
	\fontsize{10}{10}\selectfont
    \caption{Definitions of the parameters, indices, and variables used in the linear program.}
   \label{tab:model_variables}
        \centering 
   \begin{tabular}{c l } % Column formatting, 
      \hline 
      Parameters     & Description \\
      \hline 
      $N$               & Total number of vehicles designs \\
      $I$               & Total number of network nodes \\
      $P_I$             & Number of integer commodity types \\
      $P_F$             & Number of float commodity types \\
      $P = P_I + P_F$   & Total number of commodity types \\
      $T_{LP}$          & Total number of time steps in the linear program \\
      $C_{n, j, p, t}$  & ost of launching commodity type $p$, carried by vehicle $n$,\\ 
                               & to node $j$ at time $t$. \\
      $d_{n, i, p, t} $       & Demand matrix defining the supply (positive value) or demand \\
                            & (negative value) of commodity type $p$, at node $i$ at time $t$.\\ 
      $Z_{n,i,j}$           & Propellant mass fraction associated with vehicle $n$ \\
                            & travelling from node $i$ to node $j$.\\
      $c$                   & Crew consumables consumption rate \\
      $\rho$                & ISRU propellant production rate \\
      $\mu $                & ISRU maintenance supply requirement \\
      $\tau_{i,j,t}$        & Real time of flight between node $i$ and node $j$ at discrete time index $t$ \\
      $\mathcal{E}_{i,j}$   & Boolean variable defining whether the arc from node $i$ to node $j$ exists. \\
      \hline
      Index    & Description \\
      \hline 
      $n \in [0, N)$    & Vehicle  \\
      $i \in [0, I)$    & Start node  \\
      $j \in [0, I)$    & Final node  \\
      $p \in [0, P)$    & Payload type \\
      $io \in \{0,1\}$    & In-to- or out-of-arc \\
      $t \in [0,T_{LP})$& Time    \\
      \hline
      Variables    & Description \\
      \hline
      $x_{n, i, j, p, io, t}$  & Quantity of commodity type $p$, carried by vehicle $n$,\\ 
                               & from node $i$ to node $j$ at time $t$. \\
      \hline 
   \end{tabular}
\end{table} 

\begin{table}[]
    \centering
    \caption{Network nodes and arcs.}
    \begin{tabular}{c c c }
    \hline
       Node Index & Name & Arcs to \\
       \hline
        0 & Earth Surface & 0, 1 \\
        1 & Low Earth Orbit (LEO) & 0, 2 \\
        2 & Low Lunar Orbit (LLO) & 0, 1, 2, 3 \\
        3 & Lunar Surface & 2, 3 \\
        \hline
    \end{tabular}
    \label{tab:nodes}
\end{table}

\begin{table}[]
    \centering
    \caption{Payload types and their costs per unit used in the logistics formulation. Payloads 0 and 1 are integer quantities, and payloads 2 - 7 are continuous quantities.}
    \begin{tabular}{c c c }
    \hline
       Index & Payload Type & Cost, kg per unit \\
       \hline
        0 & Vehicle & $m_{\mathrm{dry}}$ \\
        1 & Crew & $m_{\mathrm{crew}}$ \\
        \hline
        2 & ISRU Plant & 1 \\
        3 & Maintenance Supplies & 1 \\
        4 & Crew Consumables & 1 \\
        5 & Miscellaneous Non-Consumable Payload & 1 \\
        6 & Oxidiser & 1 \\
        7 & Fuel & 1 \\
        \hline\
    \end{tabular}
    \label{tab:payloads}
\end{table}

\begin{figure}[]
\centering
\includegraphics[width=0.8\textwidth]{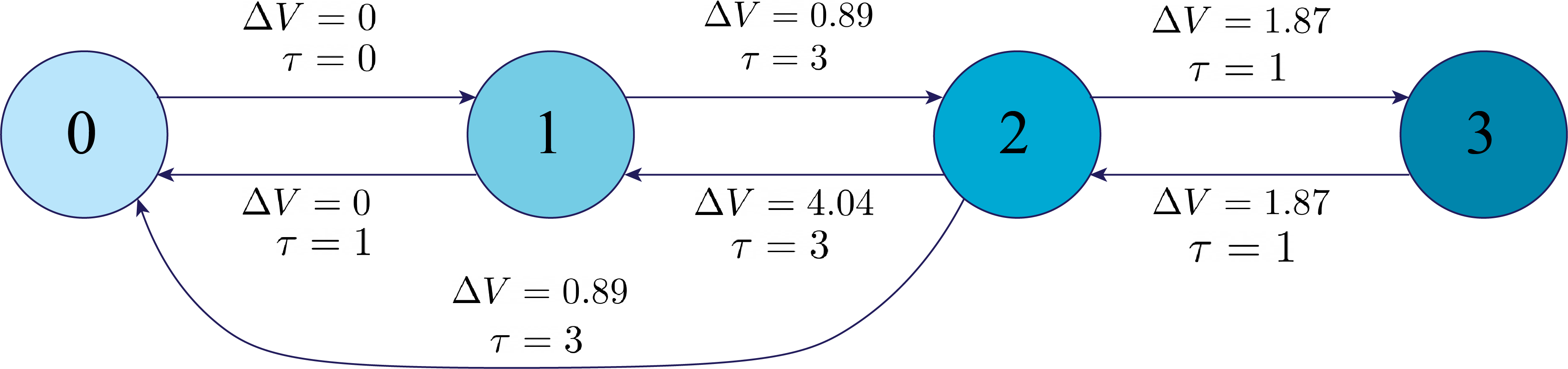}
\caption{Diagram of the network model. $\Delta V$ measured in km s$^{-1}$, values from \cite{Parker_Anderson_2014}. Lunar ascent neglects gravity losses. Time of flight measured in days. Adapted from \cite{Gollins_Isaji_Shimane_Ho_2023}.}
\label{fig:static_network}
\end{figure}

The network describes the Earth, Moon, and various orbits as ``nodes'', and the transfers between them as ``arcs''. Each arc has a cost ($\Delta V$) and a time-of-flight $\mathrm{TOF}$ associated with it. In a static network, only the nodes and the arcs representing their spatial separation are considered. The static network including the arc costs is shown in Figure \ref{fig:static_network}. The arcs from the Earth node have zero $\Delta V$ because it is assumed that a separate vehicle (the launcher) from the logistics vehicle would make these transfers. Launch vehicles are not tracked in this model, so the $\Delta V$ from node 1 to node 2 is smaller than the reverse because it is assumed that the launch vehicle performs the lunar transfer injection, with the lander performing only the lunar orbit injection in that arc. The transfer from nodes 2 to 0 has a smaller $\Delta V$ than the sum of 2 to 1 and 2 to 0 because the latter implies that the spacecraft actually inserts itself into LEO at node 1, whereas the former represents a direct atmospheric entry from the lunar return trajectory.

\begin{figure}[]
\centering
\includegraphics[width=0.7\textwidth]{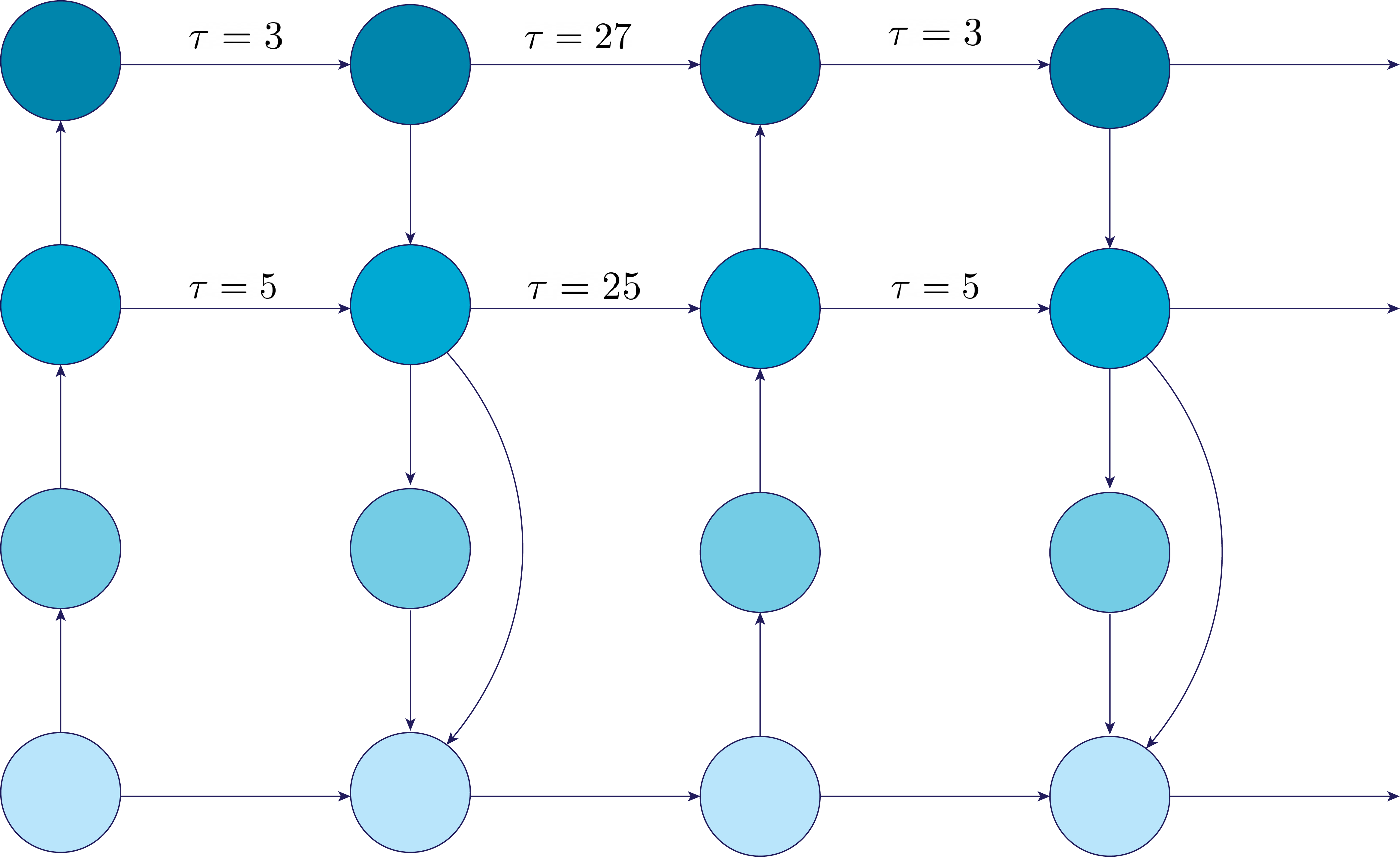}
\caption{Diagram of the time-expanded network demonstrating the direction of commodity flow on each time step, including "time of flight" $\tau$ of the holdover arcs. Time of flight measured in days.}
\label{fig:Time_Expanded} 
\end{figure}

A network becomes ``time-expanded'' by repeating the static network across many discrete time steps, with ``holdover'' arcs connecting a location to its future counterpart. This is illustrated in Figure \ref{fig:Time_Expanded}. The time expansion repeats with every second time step. That is to say, $\tau_{i,j,0} = \tau_{i,j,2}$. On these even time steps, only outbound flow is allowed. On odd time steps, return flow is allowed. A holdover arc on node 3 (lunar surface) from an even step to an odd step represents a short mission to the lunar surface, reminiscent of Apollo. The corresponding time of flight is 3 days. Meanwhile, the holdover arc in node 2 (lunar orbit) must equal the total of the time of flights for descent, surface stay, and ascent to maintain consistency, so its time of flight is 5 days. It is assumed that a long-duration mission would last some integer number of months, so the odd-to-even holdover arcs have time of flight such that 30 days have passed in total when returning to an even step. For example, a crew visiting the lunar surface takes 4 days to reach their destination. Then, a stay until the next macro-period adds 30 days. But, according to Fig. \ref{fig:Time_Expanded}, the crew must stay another semi-period (+3 days), in order to use the return-direction arcs. Then, 4 days to return to the Earth make a total of a 41-day mission. Meanwhile, a crew that stays in LLO travels 3 days to reach LLO, then stays 30 days, after which 5 days takes them to the next semi-period, followed by 3 days to return, for a consistent 41-day mission. Longer duration stays would simply add multiples of the 30-day macro-period length. 

Node 0 (Earth surface) has a holdover arc, but no associated time of flight, because time of flight is only used in consumable loss calculations and this is not considered to be relevant pre-launch.

The logistics MILP layer could be solved to optimize a variety of objectives, such as launch mass, launch cost, or fewest separate launches. For simplicity in the presented case studies, the logistics were optimized for minimized total mass to LEO, as this uses the simplest cost function. This is shown in Equation (\ref{eqn:LP_obj}). For generality, the cost coefficient $C$ is fully indexed for all nodes and all time, but with the stated objective, all coefficients are 0 unless $j=1$, and with no variation with time.

\begin{equation} \label{eqn:LP_obj}
    \min_x f(x) =  \sum_t \sum_n \sum_p \sum_j C_{n, j, p, t} x_{n, 0, 1, p, 0, t}
\end{equation}

In the following definitions, the model constraints have been arranged into 5 categories, within each of which the constraints have a shared purpose. Firstly, programmatic requirements are imposed on the commodity flow using a demand matrix $d_{n, i, p, t} $. This defines the supply or demand of commodities at specific nodes to be delivered at a specific time. Equations (\ref{eqn:LP_demand_constr_veh}) and (\ref{eqn:LP_demand_constr}) state that the difference between the total amount of commodity flowing \textit{into} a node from all others and the amount of the same commodity \textit{leaving} the node to all others, is limited by the supply (positive value $d$) or demand (negative value $d$). Note that vehicles satisfy the demand/supply of their specific type, or vehicle stacks in which they appear, $\mathcal{S}_n'$, whilst for all other payloads, only the sum of the payloads delivered by all logistics vehicles is considered. This allows commodities to transfer between vehicles but maintains proper supply rules for the vehicles themselves.

\begin{equation} \label{eqn:LP_demand_constr_veh}
    c_{1a}(x): \quad \sum_{n' \in S_n'}\sum_j \left(x_{n', i, j, p, 0, t} - x_{n', j, i, p, 1, t} \right) \leq d_{n, i, p, t} \quad \forall \; n, i, t, p = 0
\end{equation}

\begin{equation} \label{eqn:LP_demand_constr}
    c_{1b}(x): \quad \sum_n \sum_j \left(x_{n, i, j, p, 0, t} - x_{n, j, i, p, 1, t} \right) \leq \sum_n d_{n, i, p, t} \quad \forall \; i, t, p > 0
\end{equation}

The second set of constraints enforce vehicle payload capacities. Crew, which are treated as integer payloads, have a mass (cost coefficient) of $m_\mathrm{crew}$ kg per crew member, whose value is assumed to be 100 kg in this paper. Spacecraft have a cost coefficient equal to their dry mass. Payloads represented by continuous variables have a cost coefficient of 1. Note that the mass of the ISRU plants (index $p = 2$) is excluded from the capacity for holdover arcs and the lunar surface as, in ISRU-based scenarios, they are intended to remain in place on the lunar surface independent of the movement of logistics vehicles.

\begin{equation} \label{eqn:LP_pay_constr}
    c_2(x):\quad  \begin{cases}
                    m_\mathrm{crew}x_{n, i, j, 1, 0, t}+ \sum_{p=2}^5 x_{n, i, j, p, 0, t}  \leq m_{\mathrm{pay},n} x_{n, i, j, 0, 0, t} , & \forall \; n, t, (i,j): i \neq j \;  \text{or} \; i \neq 3 \\[1em]
                    m_\mathrm{crew}x_{n, i, j, 1, 0, t}+ \sum_{p=3}^5 x_{n, i, j, p, 0, t}  \leq m_{\mathrm{pay},n} x_{n, i, j, 0, 0, t} , & \forall \; n, t, (i,j): i = j = 3 \\
                \end{cases}
\end{equation}

The third set of constraints, shown in Equations (\ref{eqn:LP_ox_constr}) and (\ref{eqn:LP_fuel_constr}), impose the propellant capacity constraints.

\begin{equation} \label{eqn:LP_ox_constr}
    c_{3a}(x): \quad  x_{n, i, j, 6, 0, t}  \leq x_{n, i, j, 0, 0, t} \: \phi_n \: m_{\mathrm{prop},n} \quad \forall \; n, t, i, j 
\end{equation}

\begin{equation} \label{eqn:LP_fuel_constr}
    c_{3b}(x): \quad  x_{n, i, j, 7, 0, t}  \leq x_{n, i, j, 0, 0, t} \: (1-\phi_n) \: m_{\mathrm{prop},n} \quad \forall \; n, t, i, j 
\end{equation}

The fourth set of constraints ensure that changes in commodity quantities follow proper dynamics or conservation rules. Firstly, Equation (\ref{eq:cons_constr}) shows how crew consumables are consumed at a constant rate. In all case studies, a consumable rate of 8.655 kg per crew member per day was used, in consistency with \cite{Chen_Ho_2018}.

\begin{equation} \label{eq:cons_constr}
    c_{4a}(x): \quad x_{n, i, j, 4, 1, t} - x_{n, i, j, 4, 0, t} + c \: \tau_{i,j,t} \: x_{n, i, j, 0, 1, t} = 0 \quad \forall \; n, t, i, j
\end{equation}

Equation (\ref{eq:maint_constr}) governs the maintenance supplies associated with ISRU infrastructure. The amount of maintenance supplies per day required is proportional to the mass of ISRU infrastructure present on the lunar surface.

\begin{equation} \label{eq:maint_constr}
    c_{4b}(x): \quad \sum_n \left(x_{n, 3, 3, 3, 1, t} - x_{n, 3, 3, 3, 0, t} + \mu \: \tau_{3,3,t} \: x_{n, 3, 3, 2, 1, t} \right) = 0 \quad \forall \; t
\end{equation}

Equations (\ref{eq:ox_constr}) and (\ref{eq:fuel_constr}) describe oxidizer and fuel consumption when traveling over arcs. Holdover arcs (except the lunar surface) suffer from oxygen boil-off, which is modeled as a fractional loss-per-day $\beta$. Holdover arcs on the lunar surface ($i=3$) allow for refueling from ISRU-produced propellant, produced at a constant rate $\rho$. Transfer arcs are sufficiently short that boil-off was neglected.

\begin{equation} \label{eq:ox_constr}
    c_{4c}(x): \begin{cases}
                
                 x_{n, i, j, 6, 1, t} - (1- \beta_{n})^{\tau_{i,j,t}}x_{n, i, j, 6, 0, t} = 0  &  \forall \; n, t, (i,j):i =j \neq 3 \\[1em]

                 x_{n, i, j, 6, 1, t} - x_{n, i, j, 6, 0, t} - \rho \: \tau_{i,j,t} \sum_{n} x_{n, i, j, 2, 1, t}= 0 &  \forall \; n, t, (i,j):i =j = 3   \\[1em]

                \begin{aligned}
                 & x_{n, i, j, 6, 1, t} - x_{n, i, j, 6, 0, t} +  \\ 
                 & \phi_n Z_{n,i,j} \left(m_{\mathrm{dry},n}x_{n, i, j, 0, 0, t}  + m_{\mathrm{crew}}x_{n, i, j, 1, 0, t} + \sum_{p=2}^7  x_{n, i, j, p, 1, t}\right) = 0 
                \end{aligned} & \forall \; n, t, (i,j):i \neq j  \\[1em] 
            \end{cases}
\end{equation}

\begin{equation} \label{eq:fuel_constr}
    \begin{split}
        c_{4d}(x): \quad  & x_{n, i, j, 7, 1, t} - x_{n, i, j, 7, 0, t} + \\ 
                    & (1-\phi_n) Z_{n,i,j} \left(m_{\mathrm{dry},n}x_{n, i, j, 0, 0, t}  + m_\mathrm{crew}x_{n, i, j, 1, 0, t} + \sum_{p=2}^7x_{n, i, j, p, 1, t}\right) = 0 \quad \forall \; n, t, i, j
    \end{split}
\end{equation}

Equation (\ref{eq:other_pay_constr}) states that other commodities are simply conserved across arcs.

\begin{equation} \label{eq:other_pay_constr}
    c_{4d}(x): \quad x_{n, i, j, p, 1, t} - x_{n, i, j, p, 0, t} = 0 \quad  \forall \; n, t, i, j, p \in \{0, 1, 3, 4, 5\}
\end{equation}

The final set of constraints ensure, shown in Equation (\ref{eq:LP_zeros_constr}) that commodities only flow along arcs that exist at the current time step, and that vehicles remain within their domains.

\begin{equation}
    c_{5}(x): \quad x_{n, i, j, p, io, t} = \begin{cases}
                                            0 \quad \forall \; p, io, t,(n, i, j) \notin \mathcal{D}_n \\[1em]
                                            0 \quad \forall \; n, p, io, t, (i, j) :\mathcal{E}_{i,j} = 0 \\[1em]
                                            0 \quad \forall \; n, p, io, (i,j):i>j \; \text{and} \; t \; \text{even} \\[1em]
                                            0 \quad \forall \; n, p, io, (i,j):i<j \; \text{and} \; t \; \text{odd}  \\[1em]
                                            \text{unconstrained otherwise}
                                            \end{cases} 
    \label{eq:LP_zeros_constr}
\end{equation}

To summarize, the MILP problem solved is shown in Equation (\ref{eq:MILP_summary}).

\begin{equation}
\label{eq:MILP_summary}
    \begin{split}
        \min_x & \quad \mathrm{Eqn. \; (\ref{eqn:LP_obj})} \\
        \mathrm{s.t.} & \quad \mathrm{Eqns \; (\ref{eqn:LP_demand_constr_veh}) \sim (\ref{eq:other_pay_constr})} \\
    \end{split}
\end{equation}

The MILP-layer commodity flow optimization was carried out using the Gurobi optimization software \cite{gurobi}.

\subsection{Model Construction}
Two aspects of the commodity flow model change between iterations: the timeline, and the demand matrix. Both of these are determined by the decision vector chosen by the metaheuristic layer. Firstly, the metaheuristic chooses time stamps for each payload relative to the campaign timeline. The number of steps in the MILP timeline $T_{LP}$ is equal to two times (outbound and return) the number of unique values in the metaheuristic decision vector $\mathbf{x}(t_l)$. So, for example, in a 12-month campaign, if payload 0 launches in the third month, then $t_0 = 2$ (indexing starts at 0). The MILP, however, only cares about time steps where events happen. So if another payload 1 launches in month 7, $t_1 = 6$, then the MILP does not consider what happens on months 0, 1, 3, 4, or 5. So, the 12-month campaign timeline is mapped onto a reduced timeline, which in this case only has 4 time steps (the 2 steps where launches occur, and the half-steps allowing for return flows). Of course, the ``real'' (non-reduced) time that has passed must be tracked between the reduced steps, so that consumable calculations are made properly. So for every time step that was cut from the campaign timeline when producing the reduced timeline, 30 days are added to the ``time of flight'' of the corresponding holdover arcs.

The demand matrix is generated by creating positive demands (supplies) at the source nodes specified in the programmatic input data, and negative demands at the targets nodes, of the quantity and type of payload specified, at the time indicated in the metaheuristic decision vector. Vehicles are also supplied according to the vehicle data file. One unit of a vehicle is supplied at Earth on its first available timestep, and subsequent units are supplied with every multiple of that vehicles launch frequency. If a timestep in the reduced timeline corresponds to multiple launch frequency periods of that vehicle, then the model construction algorithm calculates how many units should be added at the next step. The MILP commodity flow model was constructed using the Pyomo python library \cite{hart2011pyomo, bynum2021pyomo}.

\begin{figure}
    \centering
    \includegraphics[width=0.9\textwidth]{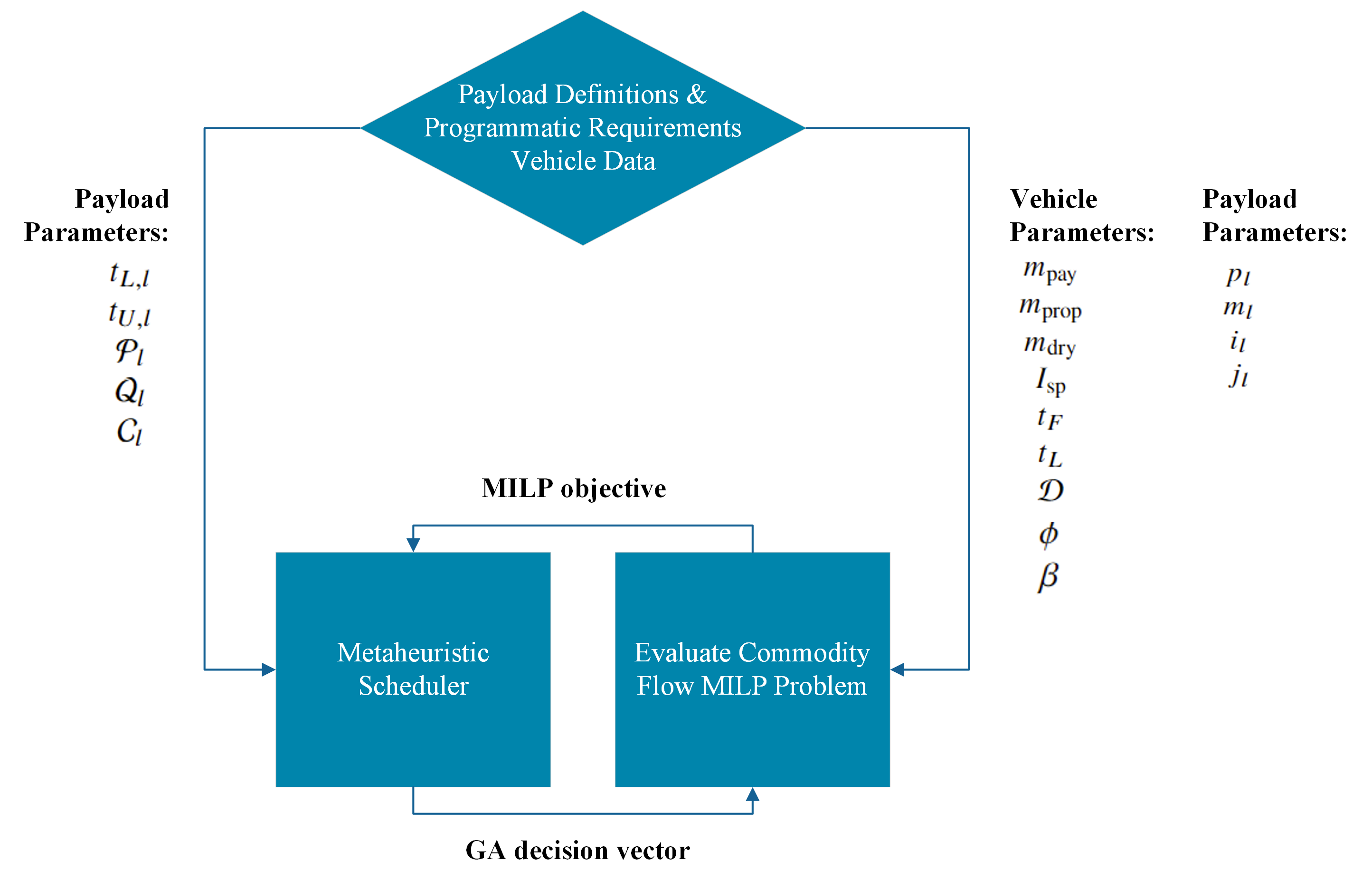}
    \caption{Flow of parameters and variables between the layers of the algorithm.}
    \label{fig:information_flow}
\end{figure}

The overall flow of information between the different layers of the algorithm is shown in Figure \ref{fig:information_flow}.

\subsection{Initial Feasible Guess Generation}
The feasible space of the problem solved by the metaheuristic layer is very sparse, so it was necessary to create an algorithm that can generate pseudo-random (although na\"{i}ve) initial feasible guesses, otherwise the feasible space is difficult to find. The guess generation method is summarised in Algorithm \ref{alg:guess_generator}. 

\begin{algorithm}
    \small
    \caption{Initial feasible guess generator pseudo-code.}
    \label{alg:guess_generator}
    \DontPrintSemicolon
    Input: Program and vehicles definitions. \\
    \textbf{x} = \textbf{0}                 \tcp*[r]{Initialize \textbf{x} with 0's.}
    \textbf{v} = rand (0, $Nv$)             \tcp*[r]{Initialize a list of vehicles assigned to each payload with random choices.}
    \For{$i \in  [0, N_P]$} { 
        \If{$\mathcal{C}_i \neq \emptyset$}{
            $\textbf{x}[i] = \textbf{x}[j] \; , \;j \in \mathcal{C}_i$
        }\Else{
            
            $t_{L,i} = \max{\left(t_{L,i}, \textbf{x}[j] \; \forall j \in \mathcal{P}_i, \textbf{x}[j]+1 \; \forall j \in \mathcal{Q}_i\right)}$ \tcp*[r]{Update lower bound $t_{L,i}$ of $\textbf{x}[i]$ to whichever is most constraining between its programmatic lower bound, or launch times of any necessary pre-cursors.}
            \If{$t_{L,i}= t_{U,i}$}{
                \textbf{x}[$i$] = $t_{L,i} = t_{U,i}$     \tcp*[r]{Only one feasible choice of launch time.}
            }\Else{
                \textbf{x}[$i$] =  rand($t_{L,i}, t_{U,i}$)     \tcp*[r]{If multiple choices are feasible, pick a random one.}
            }
            \tcp*[l]{Next, it is necessary to check that there are usable vehicles available at this time stamp.}
            \While{Valid Vehicle Check is False}{
                \tcp*[l]{Maintain a list of vehicles which are valid for this payload. Initialize with full list of vehicles:}
                Valid Vehicle List = $[0, N]$ \;                    
                \If{Payload origin node is Earth}{
                    \For{$n \in [0, N]$}{
                        \If(\tcp*[f]{Check that vehicle $n$ is available at launch time \textbf{x}[i].}){$t_{L,n} > \textbf{x}[i]$}{
                            Remove vehicle $n$ from the Valid Vehicle List.
                        }
                    }
                    \If(\tcp*[f]{If any valid vehicles remain, continue.}){Valid Vehicle List $\neq \emptyset$}{
                        \For{$n \in$ Valid Vehicle List}{
                            Payload Mass = $m_iC_i$  \tcp*[r]{Check that the payload capacity of each vehicle is not broken by adding this payload.}
                            
                            \For(\tcp*[f]{Find other payloads launching at the same time, on vehicle $n$.}){$j \in [0,i]$}{
                                \If{$\textbf{x}[i] = \textbf{x}[j]$ and $\textbf{v}[j] = n$}{
                                    Payload Mass += $m_jC_j$
                                }
                            }
                            \If(\tcp*[f]{If the payload capacity of vehicle $n$ is broken.}){Payload Mass $> m_{\mathrm{pay},n}$}{         
                                Remove vehicle $n$ from the Valid Vehicle List.
                            }
                            \If(\tcp*[f]{Check that the same vehicle is not used within a launch frequency period $t_{F,n}$}){$\textbf{v}[j] = n$ and $0 < \lvert\textbf{x}[j] - \textbf{x}[i]\rvert < t_{F,n}$}{
                                    Remove vehicle $n$ from the Valid Vehicle List.
                                    }
                        }
                    }
                }\Else{
                    \For{$n \in [0, N]$}{
                        \uIf(\tcp*[f]{If the payload is not sourced at Earth, it must be assigned to a return vehicle.}){$[3,2] \notin \mathcal{D}_n$ or $[2,0] \notin \mathcal{D}_n$}{
                            Remove vehicle $n$ from the Valid Vehicle List.
                        }   
                    }
                }
                \If{Valid Vehicle List $\neq \emptyset$}{
                    $\textbf{v}[i] =$ rand(Valid Vehicle List) \;
                    Valid Vehicle Check = True  \tcp*[r]{Successfully found a valid vehicle for payload $i$.}
                }\Else{
                    $t_{L,i} =$ rand($t_{L,i}, t_{U,i}$) \tcp*[r]{If there are no valid vehicles, try a later time step.}
                    \If{$t_{L,i} \neq t_{U,i}$}{
                        \textbf{x}[$i$] =  rand($t_{L,i}, t_{U,i}$)
                    }\Else{
                        \textbf{x}[$i$] =  $t_{L,i}$
                    }
                }
            }
        }
    }
    Output: Feasible design vector  \textbf{x}. \\
\end{algorithm}

The objective of the algorithm is to find a launch time stamp for each payload in the programmatic requirements. It starts by checking the co-payload requirements, which if present would enforce a particular time stamp. If no co-payloads are found, then it checks for precursors. The allowed launch window is then updated with the new lower bound defined by the launch time of the most-constraining precursor. Then a random time stamp is chosen from the updated launch window. Finally, the algorithm checks that a suitable vehicle is available at this time stamp. If not, the launch window lower bound is set to the previous guess, and a new random-but-later time stamp is selected. 

\section{Results}
This section will provide demonstration of both the commodity flow MILP and the metaheuristic scheduling algorithm. Then, the results of the Artemis program case studies will be presented and discussed.
\label{sec:Results}
\subsection{Logistics MILP Demonstration Case: Apollo}
\label{sec:Apollo}
Before moving on to the case studies, the commodity flow MILP is first demonstrated using a simple, single-mission Apollo study. In this simplified Apollo model, 3 astronauts are launched from Earth, with two headed for the lunar surface, and one remaining in lunar orbit. 10 kg of payload (representing a rock sample) are carried on the return leg. This programmatic data is summarized in Table \ref{tab: Apollo_payloads}, and the available vehicles are summarized in Table \ref{tab:Apollo_vehicles}. The allowed vehicle stacks are [2, 0, 1] and [0, 1]. The linear program is provided with the input decision vector $\textbf{x}(t_l) = [0, 0, 0, 0, 0]$, indicating that the mission is launched and returned within the same time step (month), corresponding to a 3 day lunar surface mission. 

\begin{figure}[]
\centering
\includegraphics[width=0.9\textwidth]{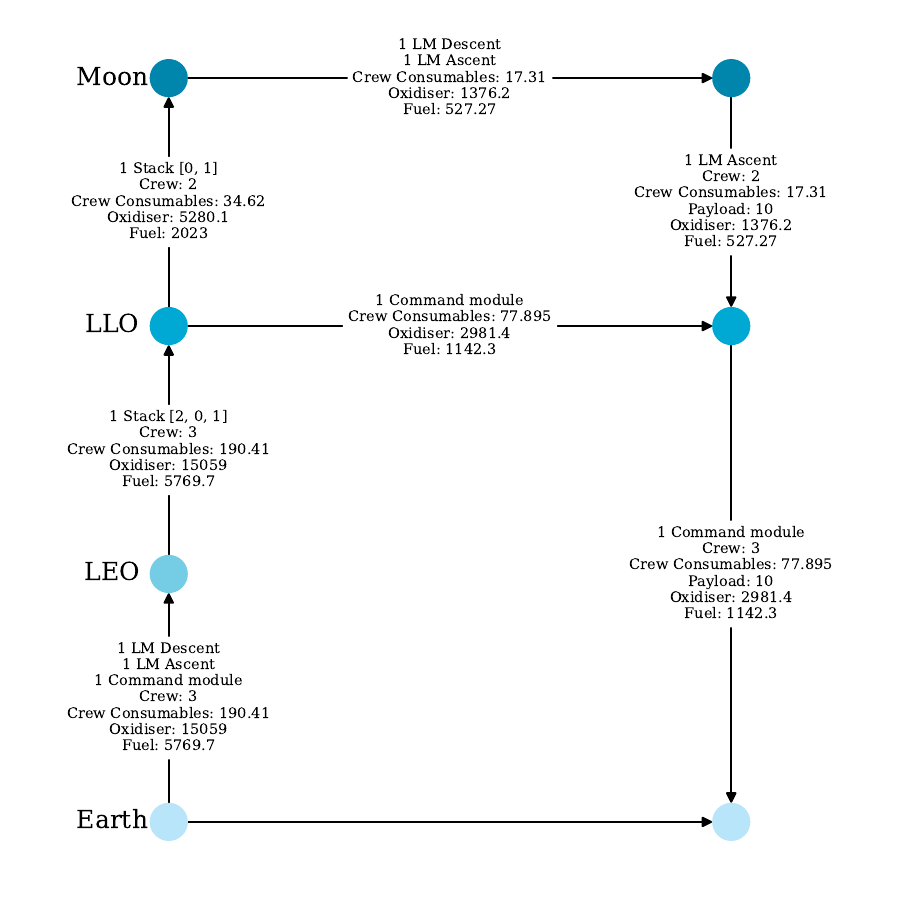}
\caption{Optimal commodity flow calculated by the MILP given the Apollo programmatic inputs.}
\label{fig:apollo_network} 
\end{figure}

The output of the linear program is summarized in Figure \ref{fig:apollo_network}. The returned objective value (total mass to LEO) is 37486 kg. This is not dissimilar to the values found by logistics MILP developed by Chen \textit{et al} \cite{Chen_Ho_2018}, with discrepancies being attributed to updates to the model and deviations in input. However, it is short of the 43572 kg mass at lunar orbit injection burn start of Apollo 11 quoted by ref. \cite{Orloff_2000}, due to the full inventory of payload and equipment not being included in this test.

\begin{table}[]
    \centering
    \caption{Apollo payload data input.}
    \begin{tabular}{l l c c c c c c c c c}
        \hline 
         $\#$ & Name & Type Index & Quantity & $i$ & $j$ & $t_L$ & $t_U$ & Co-Payloads \\
         \hline
         0 & Lunar Surface Crew & 1 & 2 & 2 & 3 & 0 & 0 & 1 \\
         1 & Lunar Orbit Crew & 1 & 3 & 0 & 2 & 0 & 0 &  \\
         2 & Surface Crew Return & 1 & 2 & 3 & 2 & 0 & 1 & 3 \\
         3 & Orbit Crew Return & 1 & 3 & 2 & 0 & 0 & 1 & \\
         4 & Sample Return & 5 & 10 & 3 & 0 & 0 & 1 & 2 \\
         \hline
    \end{tabular}
    \label{tab: Apollo_payloads}
\end{table}

\begin{table}[]
    \centering
    \caption{Apollo vehicle data input. Lunar module (L.M.) data from \cite{Isaji_Maynard_Chudoba_2020}, command module data from \cite{Orloff_2000}.}
    \begin{tabular}{l l c c c c c c c}
        \hline
         $\#$ & Name & $m_{\mathrm{pay}}$, kg & $m_{\mathrm{prop}}$, kg & $m_{\mathrm{dry}}$, kg & $I_{\mathrm{sp}}$, s & $t_F$ & $t_L$ & $\mathcal{D}$ \\
         \hline
         0 & L.M. Descent Element & 5300 & 8900 & 2217 & 311 & 1 & 0 & [0,1], [2, 2], [2, 3], [3, 3] \\
         1 & L.M. Ascent Element & 350 & 2670 & 2020 & 311 & 1 & 0 & [0,1], [3, 3], [3, 2] \\
         2 & Command Module & 22510 & 16870 & 11930 & 314.5 & 1 & 0 & [0,1], [1,1], [1, 2], [2, 2],  \\
         & & & & & & & & [2,0], [2,1], [1,0] \\
         
         \hline
    \end{tabular}
    \label{tab:Apollo_vehicles}
\end{table}

\subsection{Metaheuristic Scheduler Demonstration Case: Commercial Lunar Payloads Services (CLPS) Program}
\label{sec:CLPS}
\subsubsection{Demonstration with a Linear Cost Function} \label{sec:nonlinear}Next, the campaign scheduling aspect of the method is demonstrated by applying it to the NASA CLPS program. The CLPS program is a planned series of payloads to be delivered to the lunar surface by private companies. The payloads are a series of scientific missions that are in support of the following Artemis program. It should be noted, that the primary aim of program is to fund the development of a wide variety of commercial lunar landers for future robustness of lunar exploration campaigns, and is not optimized for the most efficient delivery of the payloads. This means, though, that it provides a relatively easy case study in which to attempt to optimize the schedule of the campaign. 

\begin{figure}[]
\centering
\includegraphics[width=1\textwidth]{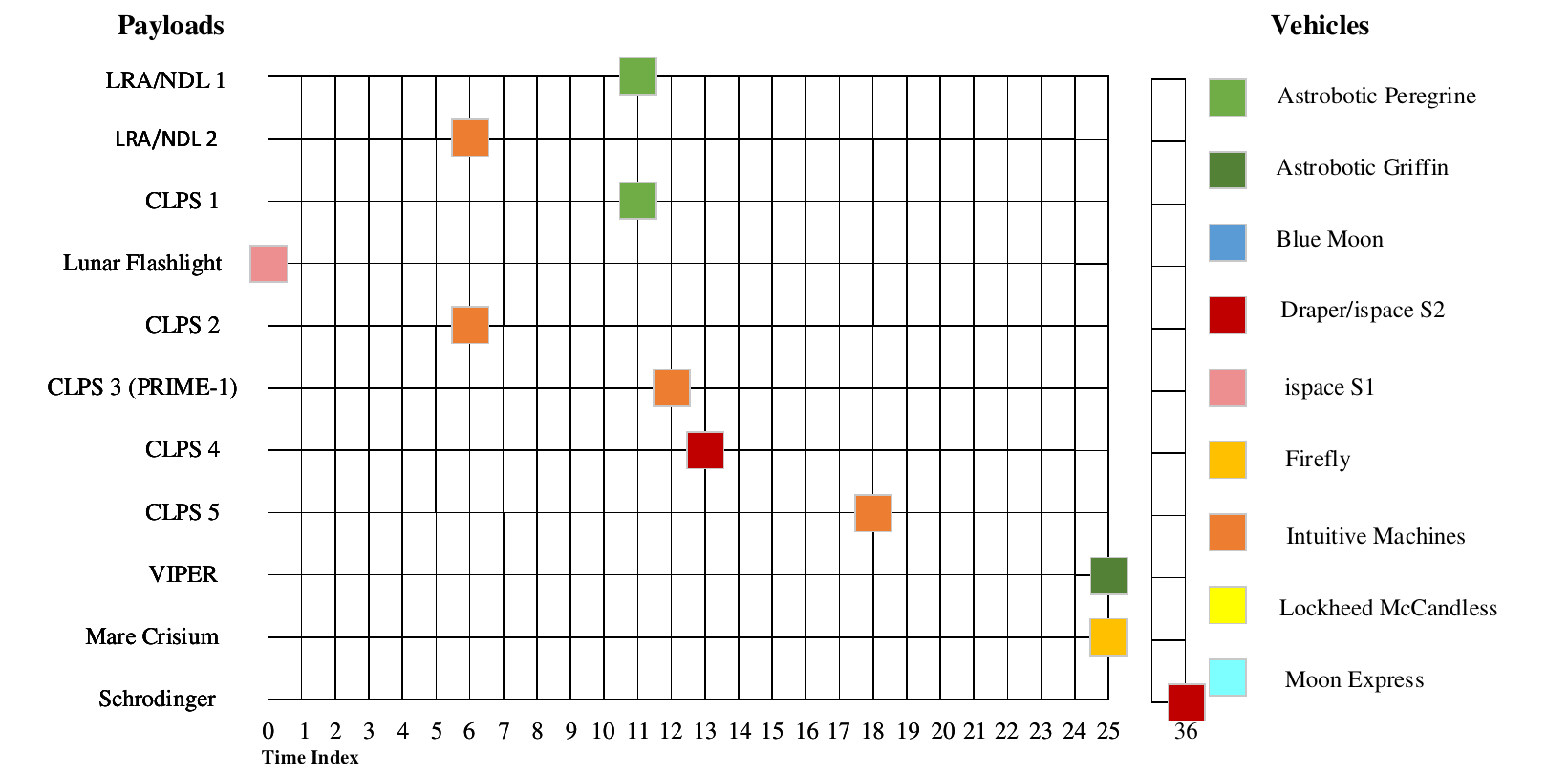}
\caption{Assumed baseline schedule and vehicle assignments of the CLPS campaign. Adapted from \cite{Gollins_Isaji_Shimane_Ho_2023}, based on assumptions made from publicly available data as of December 2022. Updated with real-world delays as of May 2023.}
\label{fig:baseline_CLPS} 
\end{figure}

In the following ConOps definitions, time $t=0$ corresponds to the month of December 2022. The list of CLPS payloads\footnote{With the addition of Lunar Flashlight. Although not a CLPS payload, it is included as it launched as a co-payload of the inaugural ispace launch.} considered in this analysis is listed in Table \ref{tab:CLPS_payloads}, and the vehicles are listed in Table \ref{tab:CLPS_vehicles}. Stacking is not allowed in this ConOps. Note that this is not an exhaustive list of CLPS vehicles, but those with sufficient publicly available data to build the model are included. Results will be compared to an assumed baseline scenario which is illustrated in Figure \ref{fig:baseline_CLPS}. 

\begin{table}[]
    \centering
    \caption{Payload data used in the CLPS analysis. Data quoted from sources where possible, and best approximations made otherwise. Time indices correspond to months with $t=0$ being December 2022. From \cite{Gollins_Isaji_Shimane_Ho_2023}. Payload mass estimations were made according to publicly available information as of December 2022.}
    \begin{tabular}{l l c c c c c c c c c}
        \hline
         $\#$ & Name & Type Index & Quantity & $i$ & $j$ & $t_L$ & $t_U$ & Co-Payloads $\mathcal{C}$ \\
         \hline
         0 & First shared payloads \cite{Warner_2020} & 5 & 14 & 0 & 3 & 0 & 12 &  \\
         1 & First shared payloads & 5 & 14 & 0 & 3 & 0 & 12 &  \\
         2 & CLPS-1 \cite{Warner_2020, Dunbar_2019} & 5 & 50 & 0 & 3 & 0 & 12 & 0 \\
         3 & Lunar Flashlight & 5 & 12 & 0 & 2 & 0 & 12 &  \\
         4 & CLPS-2 \cite{Warner_2020, Dunbar_2019} & 5 & 50 & 0 & 3 & 0 & 12 & 1 \\
         5 & CLPS-3 / PRIME-1 \cite{Vitug_2020, Dunbar_2019} & 5 & 36 & 0 & 3 & 7 & 12 & \\
         6 & CLPS-4 \cite{Dunbar_2019} & 5 & 300 & 0 & 3 & 13 & 13 & \\
         7 & CLPS-5 \cite{Dunbar_2019} & 5 & 50 & 0 & 3 & 13 & 24 & \\
         8 & VIPER \cite{Colaprete_2020, Dunbar_2019} & 5 & 430 & 0 & 3 & 12 & 25 & \\
         9 & Mare Crisium mission \cite{Potter_2021, Dunbar_2019} & 5 & 94 & 0 & 3 & 13 & 25 & \\
         10 & Schr\"{o}dinger mission \cite{Dodson_2022, Dunbar_2019} & 5 & 95 & 0 & 3 & 25 & 36 & \\
         \hline
    \end{tabular}
    \label{tab:CLPS_payloads}
\end{table}

\begin{table}[]
    \centering
    \caption{Lunar lander data used in the CLPS analysis. Data quoted from sources where possible, and best approximations made otherwise. For example, assumptions about launch mass have been made based on the assumed launch vehicle in some cases. Time indices correspond to months with $t=0$ being December 2022. Adapted from \cite{Gollins_Isaji_Shimane_Ho_2023}.}
    \scalebox{0.92}{
    \begin{tabular}{l l c c c c c c c}
        \hline
         $\#$ & Name & $m_{\mathrm{pay}}$, kg & $m_{\mathrm{prop}}$, kg & $m_{\mathrm{dry}}$, kg & $I_{\mathrm{sp}}$, s & $t_F$ & $t_L$ & $\mathcal{D}$ \\
         \hline
         0 & Astrobotic "Peregrine" \cite{Peregrine} & 90 & 720 & 470 & 340 & 12 & 1 & [0,0], [0,1], [1, 2], [2,3], [3,3] \\
         1 & Astrobotic "Griffin" \cite{Astrobotic_2021, Griffin} & 630 & 3320 & 1950 & 340 & 12 & 24 & [0,0], [0,1], [1, 2], [2,3], [3,3] \\
         2 & B.O. "Blue Moon" \cite{BlueMoon, NewGlenn} & 4500 & 6350 & 2150 & 420 & 12 & 24 & [0,0], [0,1], [1, 2], [2,3], [3,3] \\
         3 & ispace S1 \cite{ispace_2020} & 30 & 700 & 300 & 340 & 12 & 0 & [0,0], [0,1], [1, 2], [2,3], [3,3] \\
         4 & Draper/ispace S2 \cite{ispace_2021} & 500 & 3380 & 2120 & 340 & 12 & 13 & [0,0], [0,1], [1, 2], [2,3], [3,3] \\
         5 & Firefly "Blue Ghost" \cite{Firefly} & 155 & 3380 & 2470 & 340 & 12 & 13 & [0,0], [0,1], [1, 2], [2,3], [3,3] \\
         6 & I.M. "Nova-C" \cite{Berger_2021} & 100 & 1010 & 790 & 370 & 6 & 0 & [0,0], [0,1], [1, 2], [2,3], [3,3] \\
         7 & L.M. "McCandless" \cite{McCandless} & 350 & 3380 & 2270 & 340 & 12 & 48 & [0,0], [0,1], [1, 2], [2,3], [3,3] \\
         8 & M.E. MX-1 "Scout" \cite{MoonExpress_2019} & 30 & 150 & 70 & 320 & 12 & 48 & [0,0], [0,1], [1, 2], [2,3], [3,3] \\
         \hline
    \end{tabular}
    \label{tab:CLPS_vehicles}
    }
\end{table}

The schedule decision vector that corresponds to the baseline scenario of Figure \ref{fig:baseline_CLPS} is:
\begin{equation*}
    \textbf{x} = [11, 6, 11, 0, 6, 12, 13, 18, 25, 25, 36]
\end{equation*}

\begin{landscape}
\begin{figure}[]
\centering
\includegraphics[width=1.2\textwidth]{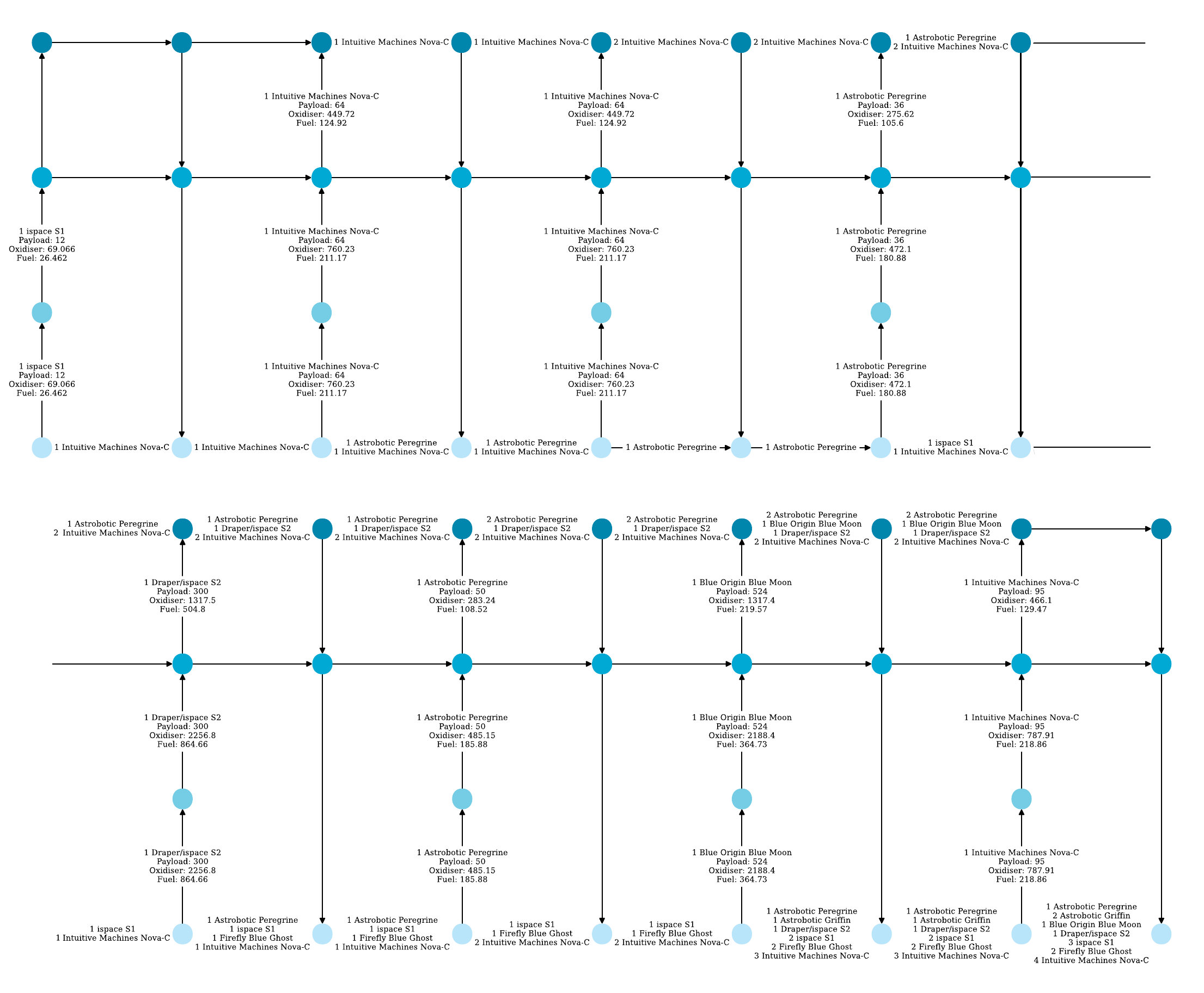}
\caption{CLPS ConOps with baseline launch schedule and optimized commodity flown.}
\label{fig:CLPS_half_optimised} 
\end{figure}
\end{landscape}

First, this vector is passed to the MILP which will solve the commodity flow optimization. The output commodity flow is shown in Figure \ref{fig:CLPS_half_optimised}, with an objective value of 19061 kg. Note that the vehicle assignments here do not match those assumed in the baseline scenario - this is because the vehicles themselves are treated as commodities, so payload-to-vehicle assignment is optimized at the MILP level.

\begin{figure}[]
\centering
\includegraphics[width=1\textwidth]{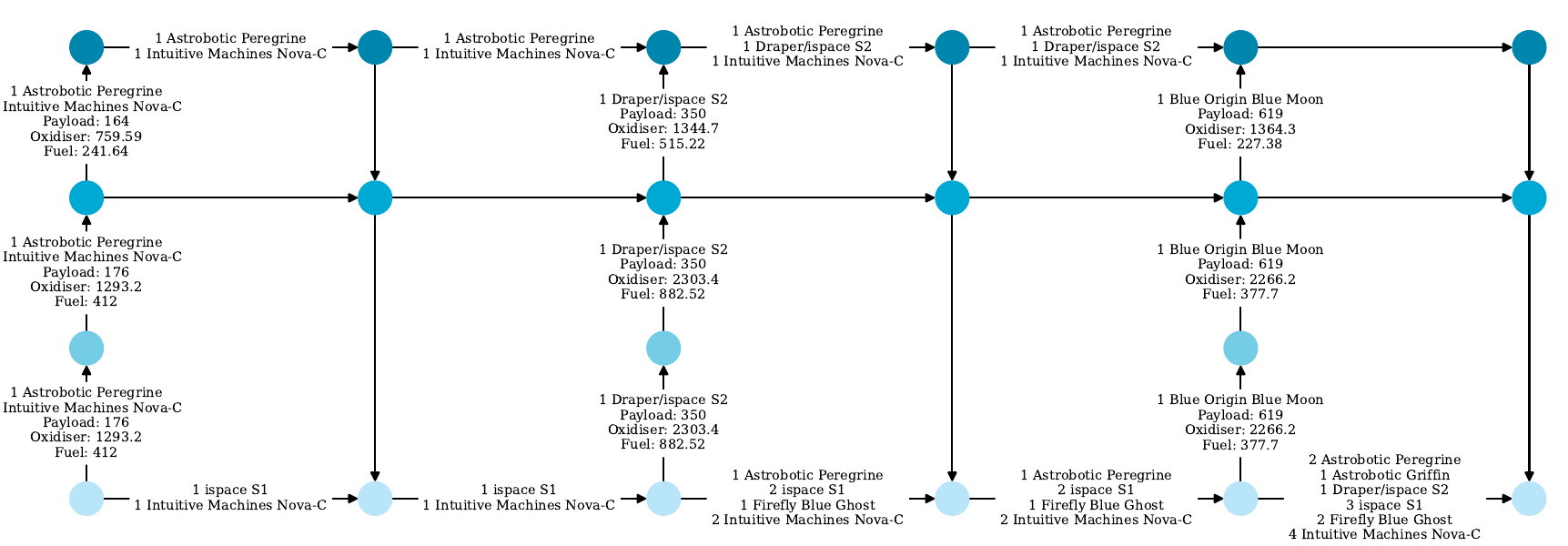}
\caption{CLPS ConOps with both optimized schedule and optimized commodity flow.}
\label{fig:CLPS_optimised} 
\end{figure}

Next, the overall campaign is optimized using the metaheuristic schedule optimizer. This was done by generating 20 initial feasible guesses using Algorithm \ref{alg:guess_generator} and using them to initialize the genetic algorithm population. 5 such populations were initialized as pygmo ``island'' objects. Islands are a framework for parallelized function evaluation with metaheuristic optimization algorithms \cite{Izzo_Ruciński_Biscani_2012}. Then, the genetic algorithm evolved the ``archipelago'' of islands through 200 generations with a mutation probability of 0.05. With each generation, the best solutions can migrate between islands, providing a form of parallelized evolution. The best-found solution through this evolution returned an objective of 14207 kg, which was found after 29 generations / 580 MILP evaluations. The commodity flow for the optimized launch schedule is shown in Figure \ref{fig:CLPS_optimised}.

As stated in the introduction, commodity flow LPs can handle scheduling with constraints no more complex than upper and lower bounds. Therefore, the CLPS case study also serves as a demonstration case for the scheduling algorithm, as the scheduling requirements are simple enough that they can be handled directly by the logistics MILP algorithm if Payloads 0 and 1 from Table \ref{tab:CLPS_payloads} are directly combined with payloads 2 and 4 such that the co-payloads requirements $\mathcal{C}$ were enforced. The new model was constructed by creating a new demand matrix in which the supply time is defined as the lower bound of the payload's availability window $t_L$ and the demand time is defined as the upper bound of the availability window $t_U$. When this new model is solved, an objective of 14207 kg is found, verifying that the scheduling algorithm found an optimal schedule.

The optimized schedule represents a mass saving of 4854 kg across the campaign. This saving is somewhat short of the equivalent of a dedicated Falcon 9 (or similar) launch. But, it would save on the cost of multiple ride-share missions of the smaller lunar logistics vehicles. This is evidenced by the fact that the total number of launches required was reduced from 7 down to 3.

\subsubsection{Demonstration with a Nonlinear Cost Function} \label{sec:nonlinear}

One of the advantages of using metaheuristic scheduler over existing MILP-based approaches is the possibility to incorporate a nonlinear cost function. To showcase this unique advantage, this subsection details a schedule optimization featuring a nonlinear cost function. As a simple example case, the same CLPS campaign was re-evaluated, this time including an exponential penalty function that penalizes against payloads launching at the same time, shown in Equation (\ref{eq:expo_penalty}). Here, the penalty function of Equation \ref{eq:meta_obj} becomes $\Gamma = \alpha \exp{(\mathcal{N}(t_l))}$, where $\mathcal{N}(t_l)$ is the number of non-unique entries listed in $\mathbf{x}(t_l)$ (i.e., penalizing launching two or more payloads at the same time step), and $k$ is a scaling coefficient. Besides the updated objective, the same constraints as previously described were maintained.

\begin{equation}
\label{eq:expo_penalty}
    \min_{t_l} \; \quad \mathcal{F}(\mathbf{x}(t_l), \mathcal{R}) + k \exp{(\mathcal{N}(t_l))}
\end{equation}

With this new objective, the optimization process of Section \ref{sec:CLPS} was repeated. The new best found decision vector was:

\begin{equation*}
    \textbf{x} = [7,  3,  7, 10,  3, 12, 13, 18, 25, 20, 35]
\end{equation*}

As expected, the algorithm returns a decision vector with minimal repeated co-manifested launches. The unpenalized commodity flow objective was 20764 kg. For the rest of the paper, the linear cost function is used as the nominal case.  

\subsection{Case Study 1: Artemis Phases 1 and 2A}
\label{sec:Art1-2A}
The first full case study expands the CLPS case to include the lunar surface missions of Artemis Phases 1 and 2A. Artemis Phase 1 is assumed to include the CLPS missions and the first crewed lunar landing (Artemis 3). Phase 2A then consists of the following surface missions up to Artemis 6. Artemis missions 3 - 5 are all 3 day missions (in terms of the timeline used in the metaheuristic scheduling algorithm, the return missions happens in the same time step as the outbound), whereas Artemis 6 is required to be at least a month long. Artemis 5 has the pressurized crew rover as a pre-requisite payload.  The list of payloads in the campaign is appended with those listed in Table \ref{tab:Art_1_2_payloads} and the list of available vehicles is appended with those in Table \ref{tab:artemis_vehicles}. The allowed vehicle stacks are [12, 10, 11], [10, 11], and [12, 11].

\begin{table}[]
    \centering
    \caption{Additonal payload data used in the Artemis phase 1 and 2A analysis. Data quoted from sources where possible, and best approximations made otherwise. Time indices correspond to months with $t=0$ being December 2022. From \cite{Gollins_Isaji_Shimane_Ho_2023}. Payload mass estimations were made according to publicly available information as of December 2022.}
    \begin{tabular}{l l c c c c c c c c c c }
        \hline
         $\#$ & Name & Type Index & Quantity & $i$ & $j$ & $t_L$ & $t_U$ & $\mathcal{P}$ & $\mathcal{Q}$ & $\mathcal{C}$ \\
         \hline
         11 & LUPEX \cite{JAXA_ISRO_2020} & 5 & 350 & 0 & 3 & 25 & 36 & & & \\
         12 & LEAP Rover \cite{Morisset_Picard_Moroso_2022} & 5 & 30 & 0 & 3 & 37 & 48 & & & \\
         13 & Unpressurised Crew Rover & 5 & 300 & 0 & 3 & 12 & 24 & & & \\
         14 & Artemis 3 Launch & 1 & 2 & 0 & 2 & 24 & 39 & & & \\
         15 & Artemis 3 Landing & 1 & 2 & 2 & 3 & 24 & 39 & & & 14 \\
         16 & Artemis 3 Ascent & 1 & 2 & 3 & 2 & 24 & 39 & & & 14 \\
         17 & Artemis 3 Return & 1 & 2 & 2 & 0 & 24 & 39 & & & 14 \\
         18 & Artemis 4 Launch & 1 & 2 & 0 & 2 & 40 & 55 & & &  \\
         19 & Artemis 4 Landing & 1 & 2 & 2 & 3 & 40 & 55 & & & 18 \\
         20 & Artemis 4 Ascent & 1 & 2 & 3 & 2 & 40 & 55 & & & 18 \\
         21 & Artemis 4 Return & 1 & 2 & 2 & 0 & 40 & 55 & & & 18 \\
         22 & ISRU demo & 2 & 300 & 0 & 3 & 37 & 49 & & & \\
         23 & Small Pressurised Rover & 5 & 4000 & 0 & 3 & 56 & 71 & & & \\
         24 & Artemis 5 Launch & 1 & 2 & 0 & 2 & 56 & 71 & 23 & &  \\
         25 & Artemis 5 Landing & 1 & 2 & 2 & 3 & 56 & 71 &  & & 24 \\
         26 & Artemis 5 Ascent & 1 & 2 & 3 & 2 & 56 & 71 &  & & 24 \\
         27 & Artemis 5 Return & 1 & 2 & 2 & 0 & 56 & 71 &  & & 24 \\
         28 & Artemis 6 Launch & 1 & 4 & 0 & 2 & 72 & 87 &  & &  \\
         29 & Artemis 6 Landing & 1 & 4 & 2 & 3 & 72 & 87 &  & & 28 \\
         30 & Artemis 6 Ascent & 1 & 4 & 3 & 2 & 72 & 87 & & 28 & \\
         31 & Artemis 6 Return & 1 & 4 & 2 & 0 & 72 & 87 &  & & 30 \\
         \hline
    \end{tabular}
    \label{tab:Art_1_2_payloads}
\end{table}

\begin{table}[]
    \centering
    \caption{Additonal lunar lander data used in the Artemis phase 1 and 2A analysis. Data quoted from sources where possible, and best approximations made otherwise. The figures regarding the Orion spacecraft differ from the real-world equivalent because this model requires that Orion deliver its crew to LLO, rather than the real-world requirement of higher lunar orbits such as near-rectilinear halo orbits. Time indices correspond to months with $t=0$ being December 2022. Adapted from \cite{Gollins_Isaji_Shimane_Ho_2023}.}
     \scalebox{0.93}{
    \begin{tabular}{l l c c c c c c c}
        \hline
         $\#$ & Name & $m_{\mathrm{pay}}$, kg & $m_{\mathrm{prop}}$, kg & $m_{\mathrm{dry}}$, kg & $I_{\mathrm{sp}}$, s & $t_F$ & $t_L$ & $\mathcal{D}$ \\
         \hline
         9 & ESA EL3 \cite{EL3_2022} & 1800 & 5580 & 2520 & 340 & 24 & 84 & [0,1], [1, 2], [2,3], [3,3] \\
         10 & ISECG lander \cite{Guidi_Haese_Landgraf_Lange_Pirrotta_Naoki_2022} & 9000 & 23660 & 9340 & 340 & 12 & 12 & [0,0], [0,1], [2,2], [2,3], [3,3] \\
         11 & ISECG ascender & 500 & 10000 & 1000 & 340 & 12 & 12 & [0,0], [0,1], [3,3], [3,2] \\ 
         12 & Orion & 11800 & 22000 & 16520 & 316 & 12 & 0 & [0,0], [0,1], [1,1], [1,2],  \\
         & & & & & & & & [2,2], [2,0], [2,1], [1,0] \\
         13 & JAXA/ISRO lander \cite{JAXA_ISRO_2020} & 350 & 3510 & 2140 & 320 & 36 & 25 & [0,1], [1, 2], [2,3], [3,3] \\
         \hline
    \end{tabular}
    }
    \label{tab:artemis_vehicles}
\end{table}

Again, an initial population of 20 random feasible guesses was generated using Algorithm \ref{alg:guess_generator}, and used to populate 5 islands. The islands were evolved through 400 generations, again with cross-migrations. In this analysis, a smaller mutation probability of 0.01 was used. This is because the decision vector is larger, due to the larger list of payloads, so there are more opportunities for mutations to occur. The smaller probability was used to balance this so that the overall mutation rate did not increase.

The full details of the best found solution are shown in Figure \ref{fig:artemis_1_2_optimised} (from here onward, unused vehicles remaining at the Earth node will not be shown for simplicity), and Table \ref{tab:Art3_6_summary} summarizes the key findings. The solution shows a quirk of the commodity flow model in that, in December 2027, the crew are launched on CLPS vehicle and rendezvous with a pre-launched Orion capsule in lunar orbit. The commodity flow solver chose this solution because it is agnostic to whether a vehicle is a crew-rated or not, and the Blue Origin lander has the highest $I_{\mathrm{sp}}$ of the available vehicles. Though this could easily be implemented in future models if necessary, it does of course indicate the value of high-$I_{\mathrm{sp}}$ crewed vehicles.

\begin{landscape}
    \begin{figure}[]
        \centering
        \includegraphics[width=1.35\textheight]{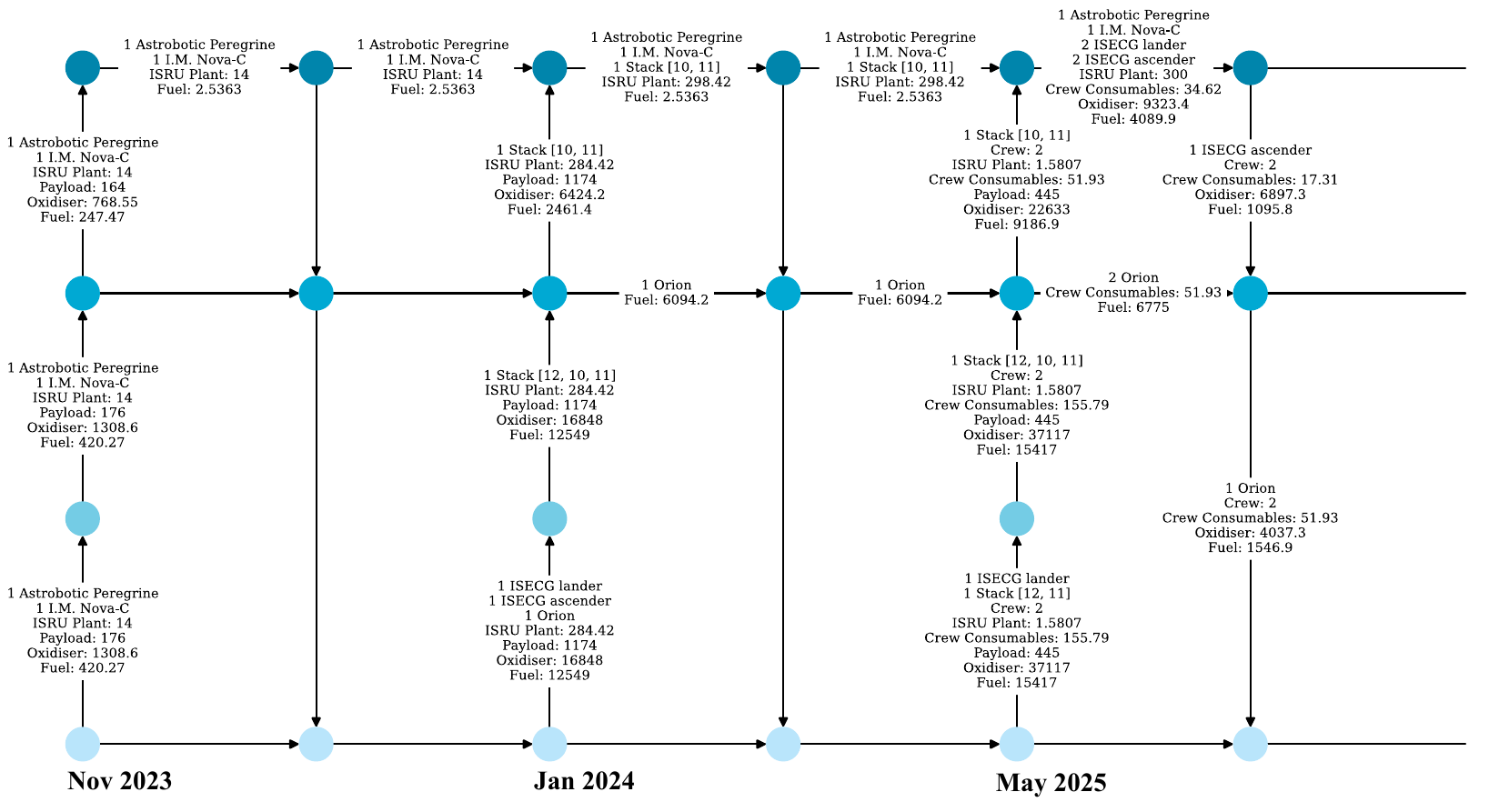}
        \caption{(continued on next page).}
    \end{figure}
\end{landscape}
\begin{landscape}
    \begin{figure}[]\ContinuedFloat
        \centering
        \includegraphics[width=1.35\textheight]{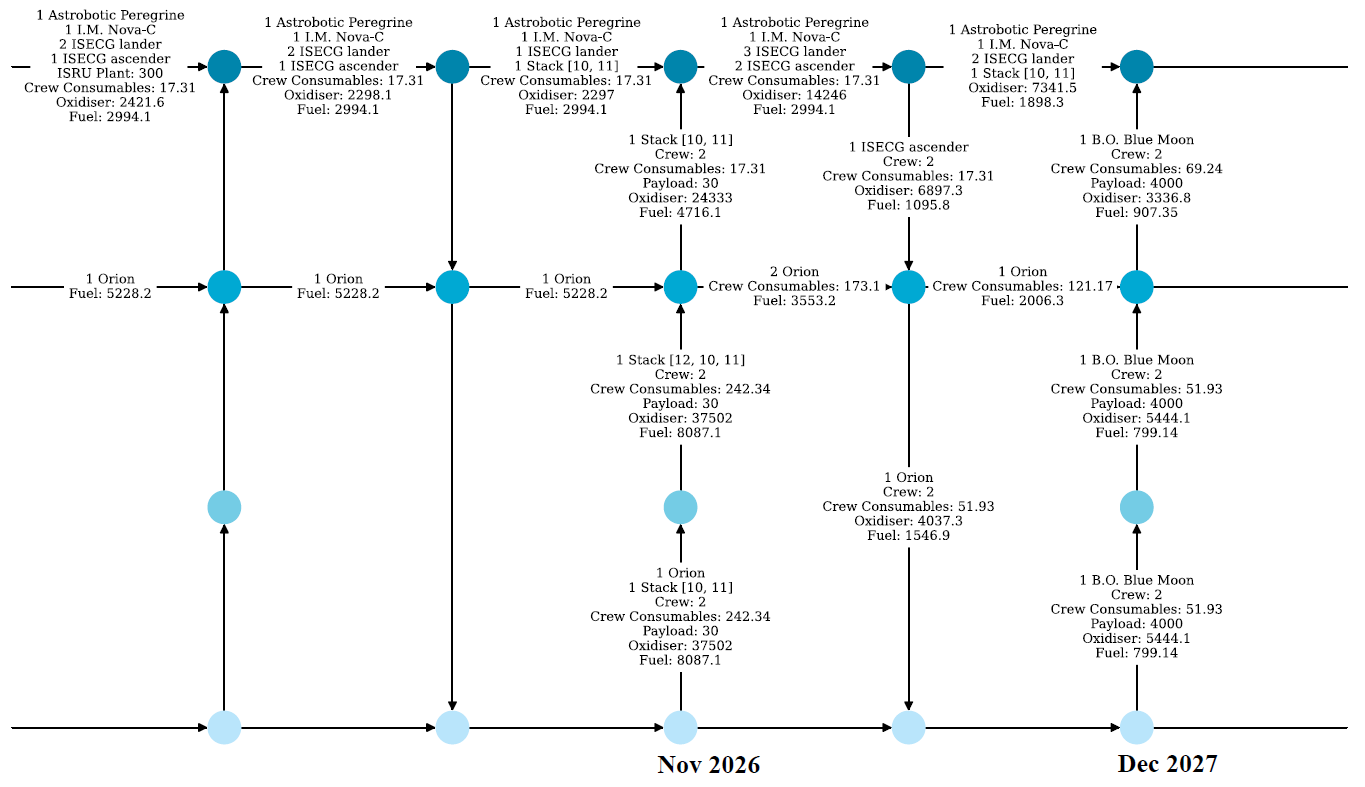}
        \caption{(continued on next page).}
    \end{figure}
\end{landscape}
    \begin{figure}[]\ContinuedFloat
        \centering
        \includegraphics[width=1\textwidth]{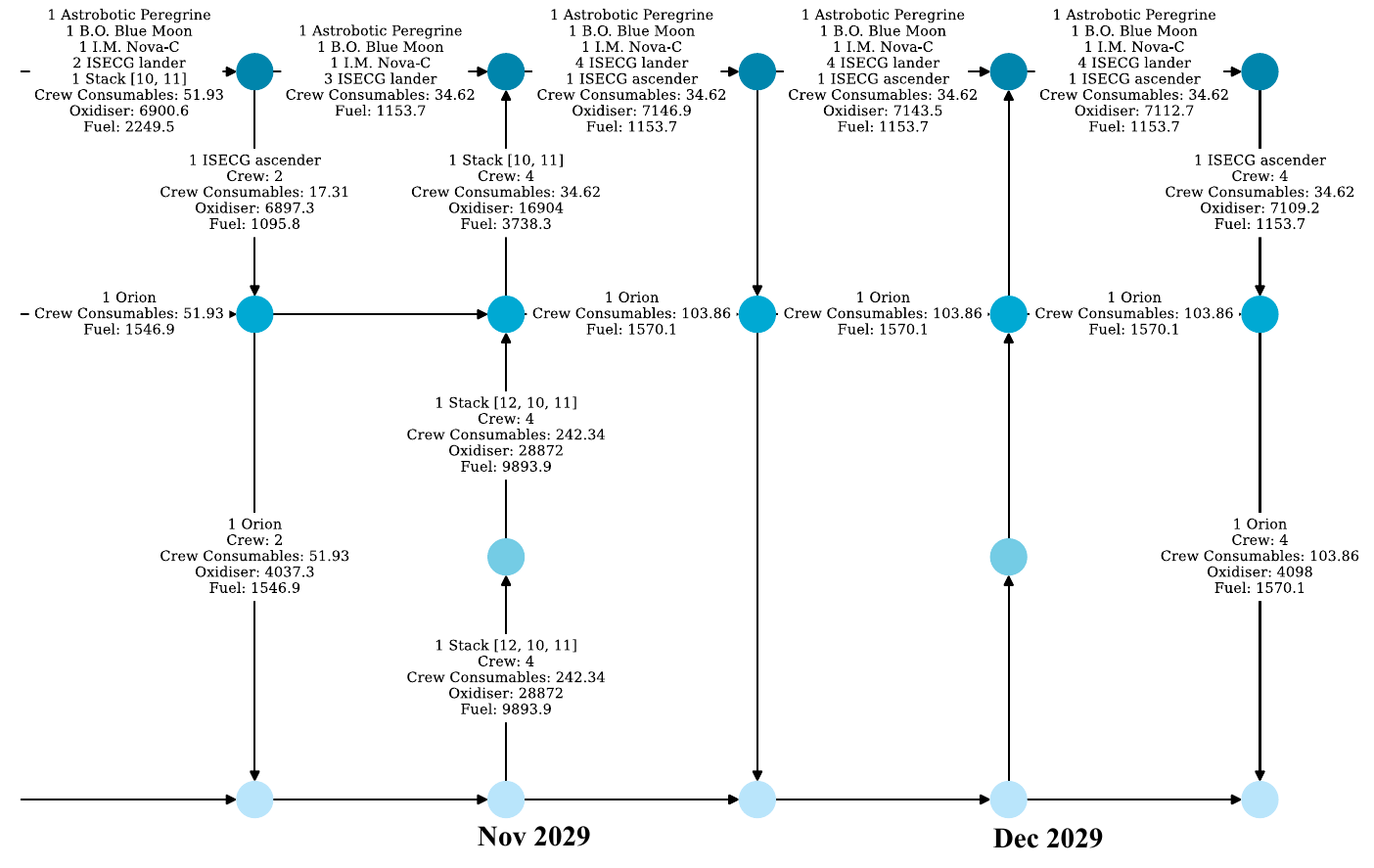}
        \caption{Artemis phase 1 and 2A ConOps with optimized schedule and optimized commodity flow.}
        \label{fig:artemis_1_2_optimised} 
    \end{figure}

\begin{table}[]
\centering
\caption{Summary of the Artemis 3 - 6 mission schedules, including pre-launched and co-manifested payloads.}
\label{tab:Art3_6_summary}
\begin{tabular}{llcccc}
\hline
\multicolumn{2}{c}{\multirow{2}{*}{\begin{tabular}[c]{@{}c@{}}Artemis \\ Mission \#\end{tabular}}} & \multirow{2}{*}{\begin{tabular}[c]{@{}c@{}}Launch \\ Date\end{tabular}} & \multicolumn{2}{l}{Pre-launched Supporting Payloads} & \multirow{2}{*}{\begin{tabular}[c]{@{}c@{}}Co-manifested \\ Supporting Payloads\end{tabular}} \\ 

\multicolumn{2}{l}{}  &   & Name  & Date   &  \\ \hline

\multirow{2}{*}{3} & Outbound  & \multirow{2}{*}{May 2025} & \multirow{2}{*}{\begin{tabular}[c]{@{}c@{}}CLPS Payloads\\ Unpressurized Rover\\ ISRU Demo\end{tabular}} & \multirow{2}{*}{\begin{tabular}[c]{@{}c@{}}Nov 2023,\\ Jan 2024\end{tabular}} & \multirow{2}{*}{CLPS Payloads} \\ 

 & Return & &  & & \\ \\
 
\multirow{2}{*}{4} & Outbound & \multirow{2}{*}{Nov 2026} & \multirow{2}{*}{} & \multirow{2}{*}{} & \multirow{2}{*}{LEAP Rover}\\

 & Return & & & & \\ 
 
\multirow{2}{*}{5} & Outbound & \multirow{2}{*}{Dec 2027} & \multirow{2}{*}{} & \multirow{2}{*}{} & \multirow{2}{*}{\begin{tabular}[c]{@{}c@{}}Small Pressurized \\ Rover\end{tabular}} \\

 & Return & & & & \\
 
\multirow{2}{*}{6} & Outbound & Nov 2029 & \multirow{2}{*}{} & \multirow{2}{*}{} & \multirow{2}{*}{} \\ 
 & Return & Dec 2029 & & & \\ \hline
\end{tabular}
\end{table}

\begin{table}
\caption{Number of times each vehicle was used in the optimized Artemis phase 1 and 2A ConOps.}
\label{tab:artemis_1_2_vehicle_usage}
\centering
\begin{tabular}{l c}
    \hline
    Name & No. times used \\
    \hline
     Astrobotic Peregrine & 1 \\
     I.M Nova-C & 1 \\
     B.O. Blue Moon & 1 \\
     ISECG lander & 4 \\
     ISECG ascender & 4 \\
     Orion & 4 \\
     \hline
\end{tabular}
\end{table}

Table \ref{tab:artemis_1_2_vehicle_usage} summarizes the number of times that vehicles were used in the optimized Artemis phase 1 and 2A ConOps. It can be seen that the optimal solution features only a small number of the available vehicles, with the optimal strategy typically being to include as many payloads in the larger capacity crewed launches as the programmatic requirements allow. It follows that the best-performing solutions are sensitive to the availability of those larger logistics vehicles. Any disruption to their availability frequency, or earliest availability, could render the best solutions infeasible. Although the unused vehicles do not contribute to the optimal results within the considered constraints, they would have usefulness in terms of robustness against certain aleatory uncertainties such as payload development delays.

\subsection{Case Study 2: Artemis Phase 2B}
\label{sec:Art2B}
A major advantage of the development of time-expanded commodity flow models for space logistics \cite{Ho_deWeck_Hoffman_Shishko_2014, Chen_Ho_2018} is the ability to model the impact that \textit{in-situ} resource utilization (ISRU) infrastructure has on the overall logistics of a campaign. In summary, resources found on the lunar surface can be taken advantage of and used to produce materials that are useful to the exploration campaign, with the objective of reducing the overall mass of material required to be launched from Earth. In particular, this method studies the impact of producing propellant from water ice on the lunar surface. 

This final case study looks to demonstrate how this can improve longer-term campaign planning and scheduling, in addition to testing a more complex set of programmatic requirements with a much larger solution space. To do this, the case study will consider the Artemis Phase 2B scenario, building towards the sustainable permanent lunar outpost of Artemis Phase 3, or the "Moon Village" concepts \cite{Biesbroek_2020}. It features the delivery of large habitation and power generation elements, and pressurized rovers. Phase 2B includes crewed missions 7 through 14, each with large sample return masses. Finally, because the impact of ISRU on the campaign was to be studied, "Mk 2" versions of the ISECG vehicles were added to the scenario, each assumed to have the same dry masses and propellant capacities but with upgraded $I_{\mathrm{sp}}$ of 370s, indicating the upgrade to ISRU-compatible cryogenic propellant. The payload capabilities of the lander were expanded accordingly with the new $I_{\mathrm{sp}}$.

In the following analysis, an ISRU infrastructure maintenance supply delivery rate of $10\%$ of the total infrastructure mass per year is used. The rate of production of \textit{in-situ} produced oxygen was calculated using the parametric sizing model developed by Schreiner \textit{et al} \cite{Schreiner_Sibille_Dominguez_Hoffman_2016}. Ref. \cite{Schreiner_Sibille_Dominguez_Hoffman_2016} offers values for reactor mass and power requirements as a function of production rate. Based on this, a relatively conservative production rate value of $2000$ kg oxygen per year per $500$ kg of reactor mass is chosen, with a power requirement of $10$ W per kg oxygen per year. The power requirement is then converted to a power system mass using data from the NASA Small Fission Power System study \cite{Gibson_Mason_Bowman_Poston_McClure_Creasy_Robinson_2015}, which offers expected specific power (W per kg) for multi-kW capable, lunar surface fission reactors. A specific power of $6.5$ W per kg is used here. With values for reactor mass, reactor power requirements, and power system mass, the total resource production rate per ISRU infrastructure mass can be calculated. With the quoted values, the resulting ISRU oxygen production rate is $\rho=0.00153$ kg oxygen per day per kg of infrastructure.

%\begin{equation}
%    \label{eq:ISRU}
%    \rho = \left(\frac{1}{\Phi_\mathrm{ISRU}} + \frac{\alpha}{\gamma}\right)^{-1}
%\end{equation}

The Artemis Phase 2B campaign is defined in this analysis by the list of payloads in Table \ref{tab:Art_2B_payloads}. These payloads are adapted from \cite{Landgraf_2021}. $t=0$ is undefined, assumed to be some date later than the December 2029 end date of the Artemis 1 and 2A analysis. The same list of vehicles from the previous case studies is used again, with the lower bound of availability removed, as this campaign is further in the future. The additional Mk 2 ISECG vehicles are listed in Table \ref{tab:artemis2B_vehicles}. Additional vehicle stacks are defined such that the Mk2 landers can stack together, as well as mixed stacks of Mk1 and Mk2 landers.

\begin{table}[]
    \centering
    \caption{Payload data used in the Artemis phase 2B analysis. Adapted from \cite{Landgraf_2021}.}
    \scalebox{0.92}{
    \begin{tabular}{l l c c c c c c c c c c }
        \hline
         $\#$ & Name & Type Index & Quantity & $i$ & $j$ & $t_L$ & $t_U$ & $\mathcal{P}$ & $\mathcal{Q}$ & $\mathcal{C}$ \\
         \hline
         0 & Power Plant element & 5 & 1500 & 0 & 3 & 0 & 48 & & & \\
         1 & Artemis 7 Crew & 1 & 4 & 0 & 2 & 0 & 48 & 0 & ~ & ~ \\ 
        2 & Artemis 7 Crew Landing & 1 & 4 & 2 & 3 & 0 & 48 & ~ & ~ & 1 \\ 
        3 & Artemis 7 Crew Ascent & 1 & 4 & 3 & 2 & 0 & 48 & ~ & 2 & ~ \\
        4 & Artemis 7 Crew Return & 1 & 4 & 2 & 0 & 0 & 48 & ~ & ~ & 3 \\
        5 & Sample return & 5 & 200 & 3 & 0 & 0 & 48 & ~ & ~ & ~ \\ 
        6 & Habitat & 5 & 4500 & 0 & 3 & 0 & 54 & ~ & ~ & ~ \\ 
        7 & Artemis 8 Crew & 1 & 4 & 0 & 2 & 0 & 54 & 6 & 4 & ~ \\ 
        8 & Artemis 8 Crew Landing & 1 & 4 & 2 & 3 & 0 & 54 & ~ & ~ & 7 \\ 
        9 & Artemis 8 Crew Ascent & 1 & 4 & 3 & 2 & 0 & 54 & ~ & 8 & ~ \\ 
        10 & Artemis 8 Crew Return & 1 & 4 & 2 & 0 & 0 & 54 & ~ & ~ & 9 \\ 
        11 & Sample return & 5 & 200 & 3 & 0 & 0 & 54 & ~ & ~ & ~ \\ 
        12 & Artemis 9 Crew & 1 & 4 & 0 & 2 & 12 & 60 & ~ & 10 & ~ \\ 
        13 & Artemis 9 Crew Landing & 1 & 4 & 2 & 3 & 12 & 60 & ~ & ~ & 12 \\ 
        14 & Artemis 9 Crew Ascent & 1 & 4 & 3 & 2 & 12 & 60 & ~ & 13 & ~ \\ 
        15 & Artemis 9 Crew Return & 1 & 4 & 2 & 0 & 12 & 60 & ~ & ~ & 14 \\ 
        16 & Sample return & 5 & 200 & 3 & 0 & 12 & 60 & ~ & ~ & ~ \\ 
        17 & Pressurised Rover & 5 & 4500 & 0 & 3 & 0 & 66 & ~ & ~ & ~ \\ 
        18 & Pressurised Rover & 5 & 4500 & 0 & 3 & 0 & 66 & ~ & ~ & ~ \\ 
        19 & Artemis 10 Crew & 1 & 4 & 0 & 2 & 24 & 66 & 17, 18 & 15 & ~ \\
        20 & Artemis 10 Crew Landing & 1 & 4 & 2 & 3 & 24 & 66 & ~ & ~ & 19 \\ 
        21 & Artemis 10 Crew Ascent & 1 & 4 & 3 & 2 & 24 & 66 & ~ & 20 & ~ \\ 
        22 & Artemis 10 Crew Return & 1 & 4 & 2 & 0 & 24 & 66 & ~ & ~ & 21 \\ 
        23 & Sample return & 5 & 200 & 3 & 0 & 24 & 66 & ~ & ~ & ~ \\ 
        24 & Artemis 11 Crew & 1 & 4 & 0 & 2 & 36 & 72 & ~ & 22 & ~ \\ 
        25 & Artemis 11 Crew Landing & 1 & 4 & 2 & 3 & 36 & 72 & ~ & ~ & 24 \\ 
        26 & Artemis 11 Crew Ascent & 1 & 4 & 3 & 2 & 36 & 72 & ~ & 25 & ~ \\ 
        27 & Artemis 11 Crew Return & 1 & 4 & 2 & 0 & 36 & 72 & ~ & ~ & 26 \\ 
        28 & Sample return & 5 & 200 & 3 & 0 & 36 & 72 & ~ & ~ & ~ \\ 
        29 & Artemis 12 Crew & 1 & 4 & 0 & 2 & 48 & 84 & ~ & 27 & ~ \\ 
        30 & Artemis 12 Crew Landing & 1 & 4 & 2 & 3 & 48 & 84 & ~ & ~ & 29 \\ 
        31 & Artemis 12 Crew Ascent & 1 & 4 & 3 & 2 & 48 & 84 & ~ & 30 & ~ \\ 
        32 & Artemis 12 Crew Return & 1 & 4 & 2 & 0 & 48 & 84 & ~ & ~ & 31 \\ 
        33 & Sample return & 5 & 200 & 3 & 0 & 48 & 84 & ~ & ~ & ~ \\ 
        34 & Fission Power Plant & 5 & 4500 & 0 & 3 & 48 & 84 & ~ & ~ & ~ \\ 
        35 & Habitat & 5 & 4500 & 0 & 3 & 48 & 84 & ~ & ~ & ~ \\ 
        36 & Artemis 13 Crew & 1 & 4 & 0 & 2 & 60 & 96 & 34, 35 & 32 & ~ \\ 
        37 & Artemis 13 Crew Landing & 1 & 4 & 2 & 3 & 60 & 96 & ~ & ~ & 36 \\ 
        38 & Artemis 13 Crew Ascent & 1 & 4 & 3 & 2 & 60 & 96 & ~ & 37 & ~ \\ 
        39 & Artemis 13 Crew Return & 1 & 4 & 2 & 0 & 60 & 96 & ~ & ~ & 38 \\ 
        40 & Sample return & 5 & 200 & 3 & 0 & 60 & 96 & ~ & ~ & ~ \\
        41 & Artemis 14 Crew & 1 & 4 & 0 & 2 & 72 & 96 & ~ & 39 & ~ \\ 
        42 & Artemis 14 Crew Landing & 1 & 4 & 2 & 3 & 72 & 96 & ~ & ~ & 41 \\
        43 & Artemis 14 Crew Ascent & 1 & 4 & 3 & 2 & 72 & 96 & ~ & 42 & ~ \\
        43 & Artemis 14 Crew Return & 1 & 4 & 2 & 0 & 72 & 96 & ~ & ~ & 43 \\ 
        44 & Sample return & 5 & 200 & 3 & 0 & 72 & 96 \\ 
         
         \hline
    \end{tabular}
    }
    \label{tab:Art_2B_payloads}
\end{table}

\begin{table}[]
    \centering
    \caption{Mk 2 ISECG vehicles used in the Artemis 2B analysis.}
    \scalebox{0.93}{
    \begin{tabular}{l l c c c c c c c}
        \hline
         $\#$ & Name & $m_{\mathrm{pay}}$, kg & $m_{\mathrm{prop}}$, kg & $m_{\mathrm{dry}}$, kg & $I_{\mathrm{sp}}$, s & $t_F$ & $t_L$ & $\mathcal{D}$ \\
         \hline
            14 & MK2 ISECG lander & 11390 & 23660 & 9340 & 370 & 12 & 0 &  [0, 0], [0,1], [2, 2], [2, 3], [3, 3] \\ 
            15 & MK2 ISECG ascender & 500 & 10000 & 1000 & 370 & 12 & 0 &  [0, 0], [0, 1], [3, 3], [3, 2] \\ 
         \hline
    \end{tabular}
    }
    \label{tab:artemis2B_vehicles}
\end{table}

20 initial feasible guesses were generated using Algorithm \ref{alg:guess_generator}. These guesses were used to initialize two metaheuristic algorithm populations. Only two population islands were used in this analysis because of the larger computing power requirements of the increased model size. A mutation probability of 0.01 was used.

\begin{figure}
    \centering
    \includegraphics[width=0.9\textwidth]{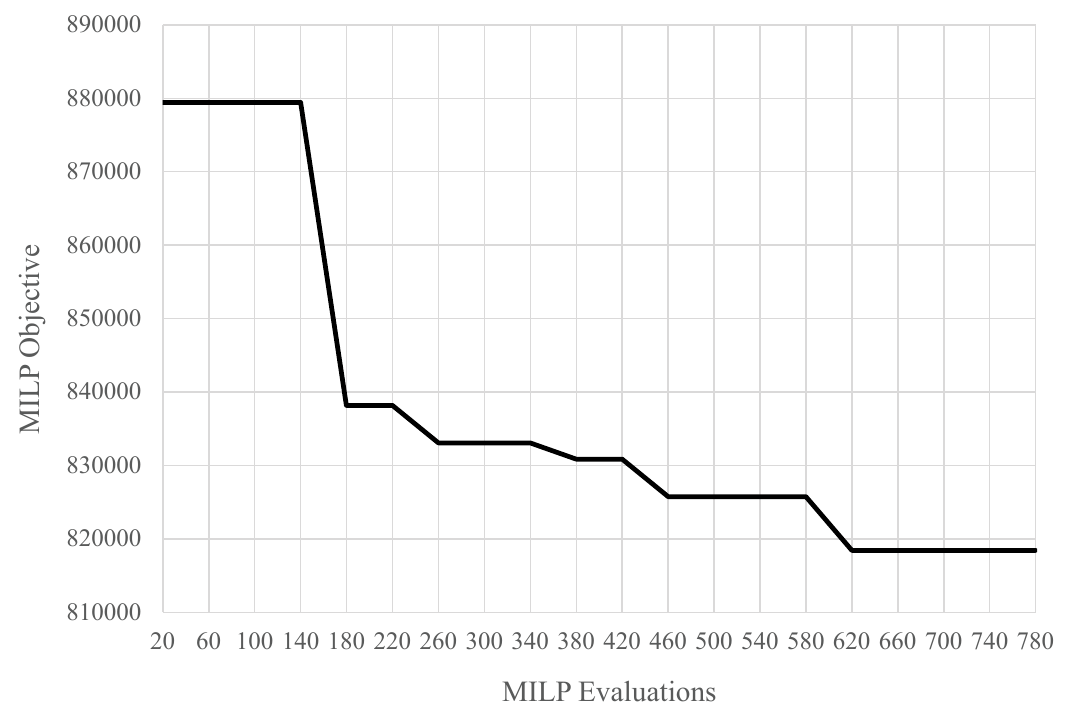}
    \caption{Evolution of the objective of the Artemis Phase 2B model.}
    \label{fig:art2B_evolution}
\end{figure}

After evolution through 40 generations, the objective value improved from 886750 kg to 819650 kg. The full commodity flow produced by this improved solution is shown in Figure \ref{fig:artemis_2B_optimised}, and the key findings are summarized in Table \ref{tab:Art_7_14_summary}. The relative saving of 67100 kg is the equivalent of $\approx 1.75$ crewed SLS Block 1B launches \cite{SLS}. The evolution that the best-found objective takes over the course of the metaheuristic optimization, shown in Figure \ref{fig:art2B_evolution}, gives some insight into the effectiveness of intermediate points between the manually-found solution and the best-found solution. Although, due to the stochasticity of genetic algorithms, these insights are not definite statements. The initial sharp improvement in the objective shows that the manually-found solution was a relatively poor one. Then, the algorithm moves through a series of small improvements, likely indicating that a number of closely related solutions exist in the design space.

\begin{table}[]
\centering
\caption{Summary of the Artemis 7 - 14 mission schedules, including pre-launched and co-manifested payloads.}
\label{tab:Art_7_14_summary}
\begin{tabular}{lllccc}
\hline
\multicolumn{2}{c}{\multirow{2}{*}{\begin{tabular}[c]{@{}c@{}}Artemis \\ Mission \#\end{tabular}}} & \multicolumn{1}{c}{\multirow{2}{*}{\begin{tabular}[c]{@{}c@{}}Launch \\ Time\end{tabular}}} & \multicolumn{2}{c}{Pre-launched Supporting Payloads}  & \multirow{2}{*}{\begin{tabular}[c]{@{}c@{}}Co-manifested \\ Supporting Payloads\end{tabular}} \\

\multicolumn{2}{c}{} & \multicolumn{1}{c}{} & Name  & Time  &  \\

\hline

\multirow{2}{*}{7} & Outbound & 46 & \multirow{2}{*}{\begin{tabular}[c]{@{}c@{}}Power Plant \\ ISRU Plant\\ Maintenance Supplies\\ Pressurized Rover\end{tabular}} & \multirow{2}{*}{\begin{tabular}[c]{@{}c@{}}17\\ \\ \\ 39\end{tabular}} & \multirow{2}{*}{\begin{tabular}[c]{@{}c@{}}ISRU Plant\\ Maintenance Supplies\end{tabular}}\\
 & Return & 48   &   &    &  \\ \\ \\\ \\ 
 
\multirow{2}{*}{8} & Outbound & 52 & \multirow{2}{*}{Habitat} & \multirow{2}{*}{48} & \multirow{2}{*}{\begin{tabular}[c]{@{}c@{}}ISRU Plant\\ Habitat\end{tabular}} \\ 
 
 & Return & 53 & & & \\
 
\multirow{2}{*}{9} & Outbound & 54 & \multirow{2}{*}{Pressurized Rover} & \multirow{2}{*}{53} & \multirow{2}{*}{} \\

 & Return & 58  &  & & \\ 
 
\multirow{2}{*}{10} & Outbound & 61 & \multirow{2}{*}{} & \multirow{2}{*}{} & \multirow{2}{*}{Fission Power Plant} \\
 
 & Return & 64 & & & \\
 
\multirow{2}{*}{11} & Outbound & 68 & \multirow{2}{*}{} & \multirow{2}{*}{} & \multirow{2}{*}{} \\

  & Return & 71 & & & \\
  
\multirow{2}{*}{12} & Outbound & 78 & \multirow{2}{*}{} & \multirow{2}{*}{} & \multirow{2}{*}{} \\

 & Return & 83 & & & \\
 
\multirow{2}{*}{13} & Outbound & 88 & \multirow{2}{*}{} & \multirow{2}{*}{} & \multirow{2}{*}{} \\

 & Return & 93 & & & \\
 
\multirow{2}{*}{14} & Outbound & 94 & \multirow{2}{*}{} & \multirow{2}{*}{} & \multirow{2}{*}{} \\

 & Return & 95 & & &  \\
 \hline
\end{tabular}
\end{table}

To analyze the structure of the problem further, it is clear that families of solutions could exist. Given a particular solution, if the groupings of payloads that launch together and the launch sequence are both maintained, it is possible to adjust the time indices of individual launches within their bounds and maintain a feasible solution, as long as vehicle supply constraints are respected. 
In this case study, the flexibility and robustness of the solution are not analyzed further and are left for future work. However, the method presented here is useful for finding good baseline solutions to the launch sequencing and payload groupings, after which minor tweaks to the schedule could be made by the program architect.

\begin{landscape}
\begin{figure}[]
\centering
\includegraphics[width=1.35\textheight]{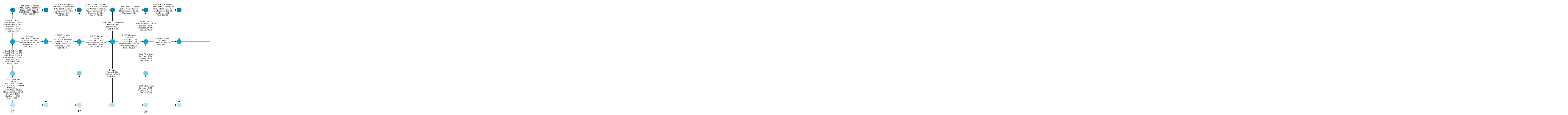}
\caption{(continued on next page).}
\end{figure}
\end{landscape}

\begin{landscape}
\begin{figure}[]\ContinuedFloat
\centering
\includegraphics[width=1.35\textheight]{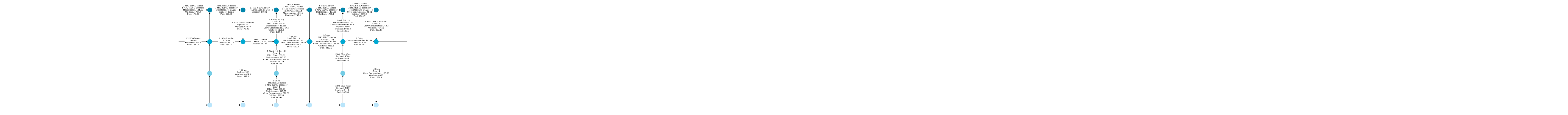}
\caption{(continued on next page).}
\end{figure}
\end{landscape}

\begin{landscape}
\begin{figure}[]\ContinuedFloat
\centering
\includegraphics[width=1.35\textheight]{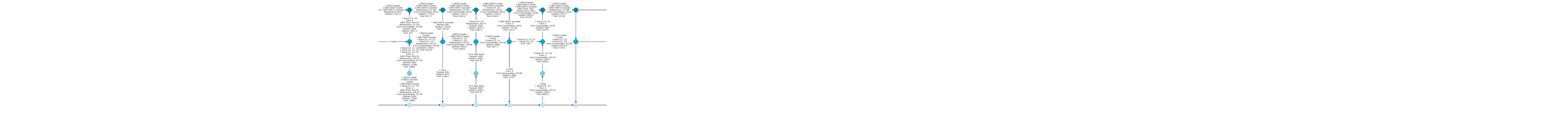}
\caption{(continued on next page).}
\end{figure}
\end{landscape}

\begin{landscape}
\begin{figure}[]\ContinuedFloat
\centering
\includegraphics[width=1.35\textheight]{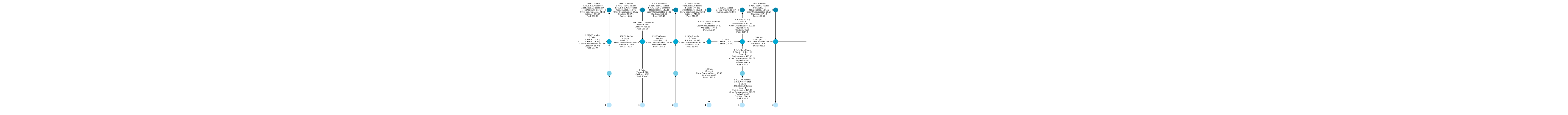}
\caption{(continued on next page).}
\end{figure}
\end{landscape}

\begin{landscape}
\begin{figure}[]\ContinuedFloat
\centering
\includegraphics[width=1.35\textheight]{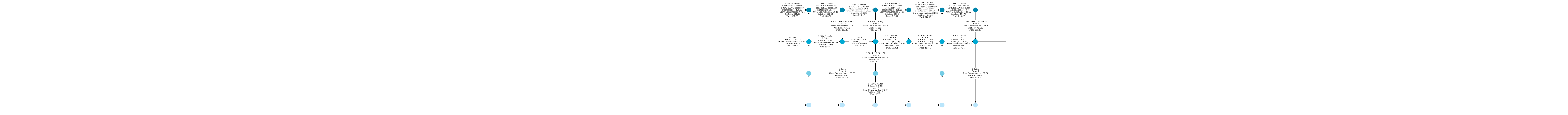}
\caption{(continued on next page).}
\end{figure}
\end{landscape}

\begin{landscape}
\begin{figure}[]\ContinuedFloat
\centering
\includegraphics[width=1.35\textheight]{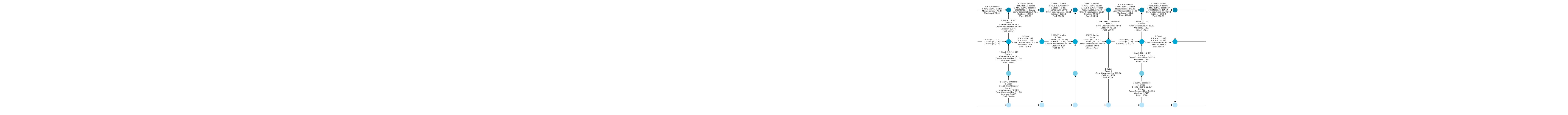}
\caption{(continued on next page).}
\end{figure}
\end{landscape}

\begin{landscape}
\begin{figure}[]\ContinuedFloat
\centering
\includegraphics[width=1.35\textheight]{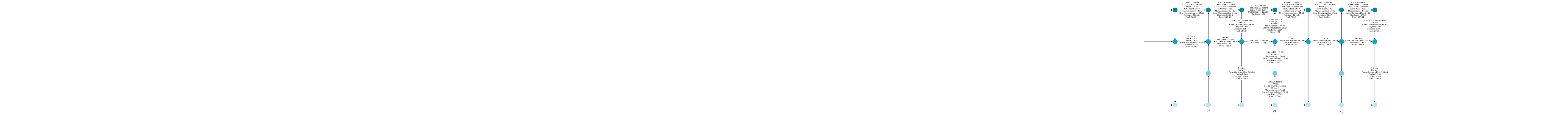}
\caption{Artemis phase 2B ConOps with improved schedule and optimized commodity flow.}
\label{fig:artemis_2B_optimised} 
\end{figure}
\end{landscape}

As expected, it can be seen that the commodity flow optimizer ships ISRU infrastructure and maintenance supplies with the first mission of the campaign. In this particular solution, ISRU infrastructure is packed into the excess payload capability of logistics vehicles that are traveling with campaign payloads already. Therefore, the exact amount of ISRU infrastructure delivered to the lunar surface is dependent on the payload capability of the vehicles. Uncertainties in the logistics vehicle performance, such as engine $I_{\mathrm{sp}}$ or structure mass ratio, would propagate through to variation in available ISRU infrastructure, and therefore creating uncertainty in surface propellant availability. Further maintenance supplies are shipped with later missions. The inclusion of ISRU infrastructure prompts the scheduler to place the later missions later in their allowed windows, so that more propellant can be produced \textit{in-situ} and reduce the amount required to be launched from Earth. For comparison, with ISRU capability disabled, the same solution produces an objective value of 845510 kg, demonstrating the net impact of ISRU on the logistics of the Artemis 2B scenario.

\section{Conclusion}
This paper has presented a method for finding feasible exploration campaign plans subject to programmatic requirements, and in particular finding schedules that produce optimal commodity flow such that total launch mass across the entire campaign is minimized. By using a mixed-integer linear program to solve the optimal commodity flow for a given campaign plan, a genetic algorithm can then be used to find improved schedules with potentially nonlinear weighting. The scheduling algorithm was demonstrated with various example cases. In the case studies, it was found that medium-to-large class lunar landers are most effective in optimizing lunar exploration logistics and, in cases of low availability of those landers, smaller and readily available landers provide useful, though sub-optimal, alternatives. 

\section*{Appendix A: MILP Solver Performance}
\begin{figure}
    \centering
    \includegraphics[width=0.9\textwidth]{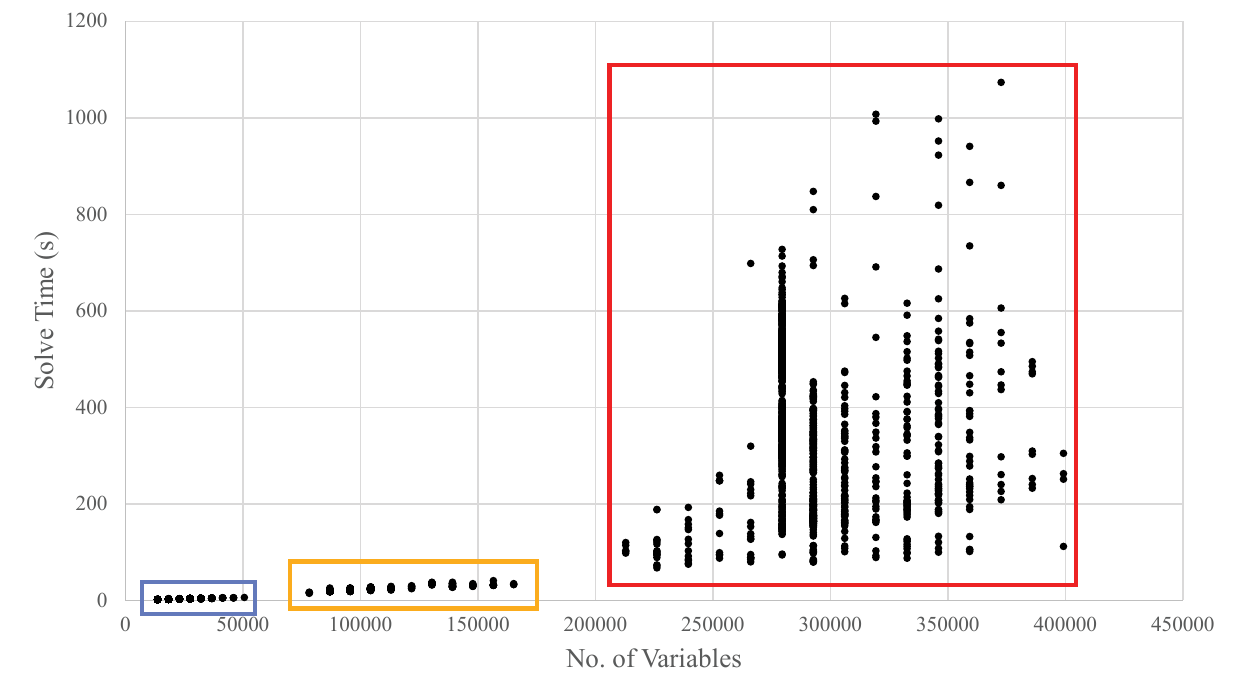}
    \caption{Solve times of the MILP versus the model size, measured across the CLPS (blue), Artemis Phase 1 - 2A (orange), and Artemis Phase 2B (red) case studies.}
    \label{fig:MILP_perf}
\end{figure}

The commodity flow problem solve time was measured for each function evaluation of the case studies discussed in this paper. Figure \ref{fig:MILP_perf} plots the solve time versus the number of variables present in each model. 5 islands in parallel of CLPS and Artemis Phase 1 - 2A models were solved using a PC with an 8-core Intel Core i7-10700 2.90GHz CPU  and 32 GB RAM, and 2 islands in parallel of Artemis 2B models were solved using the Georgia Tech PACE cluster \cite{PACE} with a 24-core Dual Intel Xeon Gold 6226 2.7 GHz CPU node and 32 GB RAM. The average solve time and the variance in the solve time increase exponentially with the model size. An all-in-MILP approach to the campaign scheduling problem would require a factor of $\approx P\times T$ more variables than those the models tested in Figure \ref{fig:MILP_perf}, because it becomes necessary to index every payload in the campaign separately in order to formulate sequencing constraints, and it also becomes impossible to trim the MILP timeline to just the time stamps in the potential solution. The solve time for such a model would be unpredictable and likely extremely long. In addition, memory allocation issues would eventually be encountered. The method presented in this paper avoids this issue by repeatedly solving reduced commodity flow problems.

\section*{Acknowledgment}
The authors acknowledge Yuri Shimane and Masafumi Isaji for their helpful feedback and suggestions. 

\bibliography{main}

\end{document}